\input ntmacros.tex
\phantom{a}
\vskip.22in
\input ntmacros.tex
\input epsf
{{\small \centerline{{\bf {REGULAR SEIFERT SURFACES AND}}}}
{{\small \centerline{{\bf VASSILIEV KNOT INVARIANTS}}}}
\bigskip
{\baselineskip=10pt{\msmall
\centerline{{\small E}{\msmall FSTRATIA } {\small K}{\msmall ALFAGIANNI\ 
{\msmall AND}\ \  {\small X}{\msmall IAO}-{\small S}{\msmall ONG}
{\small L}{\msmall IN }}}}

\bigskip

\vskip.15in
{\baselineskip=11pt {\footnote { }{ {\msmall 

\p The first author was supported by NSF grant DMS-9626140, and the second author

\p by NSF grants DMS-9400806 and DMS-9796130.

\p 1991 Mathematics Subject Classification: 57M25.

\p This version:  November 1999.}}}}

\centerline {\B CONTENTS}}
\medskip 
\smallskip
{\msmall \p{\bf  Introduction}
\medskip
\p{\bf \S 1.} {\bf Vassiliev knot 
invariants and Gusarov's {\it n-}triviality}
\medskip
\p{\bf \S 2.} {\bf Regular Seifert surfaces of a knot}
\smallskip
{2a.} Generalities

{2b.} Good position of a band in a regular spine

{2c.}  Lower central series and curves on surfaces
 
{2d.} An example illustrating the good position of a band

\medskip
\p{\bf \S 3.} {\bf Commutator half bases and Vassiliev invariants}
\smallskip
{3a.} Definitions and the main result

{3b.} Nice arcs and simple commutators

{3c.} Bad sets and good arcs  
 
{3d.} Conflict sets and products of good arcs

{3e.} The reduction to nice arcs
\medskip
\p{\bf \S 4.} {\bf Vassiliev invariants as obstructions to
 {\it n-}sliceness}
\medskip
\p{\bf \S 5.} {\bf More types of {\it n-}trivial knots}
\smallskip
{5a.} Definitions

{5b.} Mixed type knots
\medskip
\p{\bf \S 6.} {\bf Discussion and questions}
\medskip
\p{\bf References}}
\medskip

\bigskip
\bigskip 
\centerline {\B INTRODUCTION}}
\medskip

It is well known ([Bi], [BN1])
that the Vassiliev knot invariants (or knot invariants of finite type) 
give a coherent 
and systematic
way to view all the polynomial knot invariants (quantum invariants).
Although both classes of the invariants have been studied
extensively in the recent years, the following fundamental question is open:
What geometric knot properties are Vassiliev invariants
detecting? In  particular to what extent does the vanishing of 
Vassiliev invariants approximate the unknottedness of a knot?

In the present paper we
address the first question mentioned above. 
We show that
the invariants are obstructions
to a knot's bounding a {\it regular} Seifert surface whose complement
looks, modulo the lower central series
of its fundamental group, like the complement
of a null-isotopy. See Theorems 3.2 and 5.4.

Milnor's ${\bar \mu}$-invariants  ([Mi1], [Mi2]) are integer link
concordance invariants
defined in terms of algebraic data extracted from the link group.
By [BN2], [L2] and [HM], these invariants are
the only link concordance invariants that are of finite type.
Thus they
can
be thought of as
the link homotopy 
(or link concordance)
counterpart of Vassiliev's theory. 
As shown in
[Co] and [O], Milnor's
invariants detect exactly when the complements
of links are cobordant as manifolds by a cobordism
which looks, modulo the lower central series
of its fundamental group, like the complement
of a link concordance.
Our main results here
are parallel to results of
[Co] and [O] and indicate that the two classes of invariants
can be put under a unifying geometric framework.

The geometric notion
that captures the vanishing of Milnor's
invariants is
the notion of {\it $n$-sliceness}
or {\it $n$-null-concordance} ([Co],[O]).
Before we state the main results of this paper,
let us introduce some notation and terminology. We will say that
a Seifert surface $S$ of a knot $K$ is 
{\it regular} if it has a spine $\Sigma$ whose embedding in $S^3$,
induced by the embedding $S\subset S^3$, is
isotopic to the standard embedding of a bouquet of circles. Such a spine
will be called a {\it regular spine} of $S$. 
In particular, $\pi=\pi_1(S^3 \setminus S)$
is a free group.
For each circle in a regular
spine, we may push it off the surface $S$ slightly in the normal directions. 
If these push-offs are all disjoint (this can be achieved by either pushing
circles off to different sides of $S$ or in different distances away from $S$),
we get a link $L^{\bar\epsilon}$. 
For ``appropriate" $L^{\bar\epsilon}$
(see before the statement of Theorem 4.3 for details)
we obtain that
the Vassiliev
invariants of the boundary knot
$K=\partial S$ are obstructions to the null-sliceness for the links $L^{{\bar
\epsilon}}$. This in turn gives a close relation
between the Milnor invariants of the links
 $L^{{\bar
\epsilon}}$ and the Vassiliev invariants of $K$ (Theorem 4.4).
More precisely, we show that there exists
a sequence of natural numbers $\{ l(n)\}_{n\in \N}$, with
${\displaystyle { l(n) > {\rm log}_2({{n-5}\over 144})}}$,
such that the following is true:

\smallskip
\p {\bf Theorem 1.} {\sl 
Let the notation be as above. Assume that $L^{{\bar \epsilon}}$ is $n$-slice.
Then, all the Vassiliev invariants of orders $\leq l(2n-1)$
vanish for $K$.}
\medskip
The sequence $\{ l(n)\}_{n\in \N}$ is defined in $\S 3$.
A
key idea we will introduce is to define 
{\it $n$-unknotted}
Seifert surfaces. Roughly speaking, these are surfaces
whose complement 
looks, modulo certain terms of the lower central series
of its fundamental group, ``simple"
(e.g. like the complement of a Seifert surface of a trivial knot).
A knot bounding such a surface is called
{\it $n$-unknotted}. Our main result in this paper
shows that certain Vassiliev invariants of a knot
obstruct to { $n$-unknottedness.} (For details see Theorem 5.6.)
A particularly 
interesting class of surfaces that we study, 
that is used to prove Theorem 1 above, is that of
{\it $n$-hyperbolic}.
Roughly speaking,
a regular Seifert surface $S$ is called $n$-hyperbolic, for some $n\in \N$, 
if its complement looks, modulo the first
$(n+1)$ terms of its fundamental group, like the complement of a disc.
These surfaces are introduced in $\S 3$.
We prove the following:
\smallskip
\p {\bf Theorem 2.} {\sl If $K$ is $n$-hyperbolic,
then all the Vassiliev invariants of orders
$\leq l(n)$ of $K$ vanish. In particular,
if $K$ is $n$-hyperbolic for all $n\in \N$
then all the Vassiliev invariants of $K$ vanish.} 
\medskip
To motivate our notion of {\it n-hyperbolicity},
let us recall that a link is {\it n-slice} if
its components bound disjoint smooth surfaces $V_i$
in the 4-ball $B^4$ and such that each curve on $V$,
if pushed off into the complement, lies in the $n$-th
term of the lower central series of $\pi_1(B^4\setminus\cup V_i)$.
The results in [Co] and [O] and in [L1] motivated our definition above.
Our terminology though, is motivated by the
form of the Seifert matrix of the knot. See Remark
5.3.

In [Gu], Gusarov studied knots with trivial 
invariants of bounded orders by means of combinatorics
on knot projections. His notion
of {\it n-triviality} relates a knot
with vanishing finite type invariants
to the unknot, via crossing changes on knot diagrams.
This process has to take into account all
the intermediate knots that
the knot passes through. Our notion on the other hand,
involves only the knot in question and its complement,
and it is, thus, closer to the spirit
of classical geometric topology.

An interesting problem, arising from this work,
is the question of whether the
notion of $n$-hyperbolicity provides 
a complete geometric characterization of $n$-trivial knots.
We conjecture that this is the case. More precisely
we have the following: 
\medskip
\p {\bf Conjecture 3.} {\sl 
A knot $K$ is $n$-trivial for all $n\in \N$
if and only if it is $n$-hyperbolic
for all $n\in \N$.}
\medskip
The reader is referred to $\S 6$ 
for  further discussion.

Let us now describe in more detail the contents of the paper and
some of the ideas that are involved in the proofs of the main results.

In $\S 1$ we recall basic facts about 
Vassiliev invariants and the results from [Gu]
that we use in subsequent sections. 

In $\S 2$ we study regular Seifert surfaces of knots.
We introduce the notion of {\it good position} for bands
in projections of Seifert surfaces.
Let $B$ be a band of a regular Seifert surface $S$,
which is assumed to be in band-disc form,
and let $\gamma$ denote the core of $B$.
Also, let $\gamma^{\epsilon}$
denote a push-off of $\gamma$.
The main feature of a projection of $S$ with respect
to which $B$ is in {\it good position} is the following:
We may find a word $W$, in the free generators of 
$\pi=\pi_1(S^3 \setminus S)$,
representing $\gamma^{\epsilon}$
and such that every letter in $W$ is realized by
a band crossing in the projection.

In $\S 3$ we introduce
{\it n-hyperbolic} regular Seifert surfaces and we prove 
Theorem 2.
The special projections of $\S 2$ allow us  to connect Gusarov's 
notion of {\it n-triviality} to an algebraic { n-triviality}
in $\pi$, and exhibit 
a correspondence between geometry in $S^3\setminus S$
and algebra in $\pi$. 
Let us explain this in some more detail.
By [Gu], to prove Theorem 2 it will be enough to show
that an {\it n-hyperbolic}
knot has to be { $l(n)$-trivial}.
Showing that a knot is {\it k-trivial}
amounts to showing that it can be unknotted
in $2^{k+1}-1$ ways by changing crossings in a fixed projection.
See Definition 1.2.
Having  the projections of $\S 2$ at hand,
the main step in
the proof of Theorem 2 becomes showing the following:
If $\gamma$ is the core of a band $B$ in {\it good position}
and $\gamma^{\epsilon} \in \pi^{(m+1)}$,
then we can trivialize $B$
in $2^{l(m)+1}-1$ ways (for the precise statement see Proposition 3.3).
Here $\pi^{(m+1)}$ denotes the $(m+1)$-th
term of the lower central series of $\pi$.

The proof of Proposition 3.3 is based on
a careful analysis of the geometric combinatorics of projections of
the
sub-arcs of $\gamma$ representing simple commutators.
We show that eventually $\gamma$ may be decomposed
into a disjoint union of ``nice" arcs
for which the desired conclusion
follows by Dehn's Lemma.

In $\S 4$ we start by recalling the Definition
of  {\it n-slice} and a few results from [Co] and [O].
Then we show how Theorem 1 follows from Theorem 2.

In $\S 5$ we introduce $n$-unknotted
Seifert surfaces and we prove 
various generalizations
Theorem 2. See Theorems 5.4 and 5.6.

Finally, in $\S 6$
we discuss some consequences and state
some questions that arise from this work.
\smallskip
\p {\it Acknowledgment.} We thank James Conant
for a careful reading 
of this paper. His examples (see Remarks 3.3* and  3.14*) 
were crucial to  us in clarifying the subtleties of the layouts of commutators. 
Also, we thank the referee for his/her comments that were very helpful
in improving the exposition of the manuscript.
\bigskip
\medskip
\centerline {1.{\B \uppercase { Vassiliev knot invariants and Gusarov's}} 
{\it n-}{\B \uppercase { triviality}}}
\medskip
A {\it singular knot} $K\subset S^3$ is an immersed curve
whose only singularities are finitely many transverse double points.

Let ${\cal K}_n$ be the rational vector space generated by the set of ambient 
isotopy
classes of oriented, singular knots with  exactly $n$ double points.
In particular ${\cal K}={\cal K}_0$ is  the space generated by
the set of isotopy
classes of oriented knots.

A knot invariant  $V$ can be extended to an invariant of singular knots
by defining

\medskip
\centerline{\epsfysize=.4in\epsfbox[49 56 566 735]{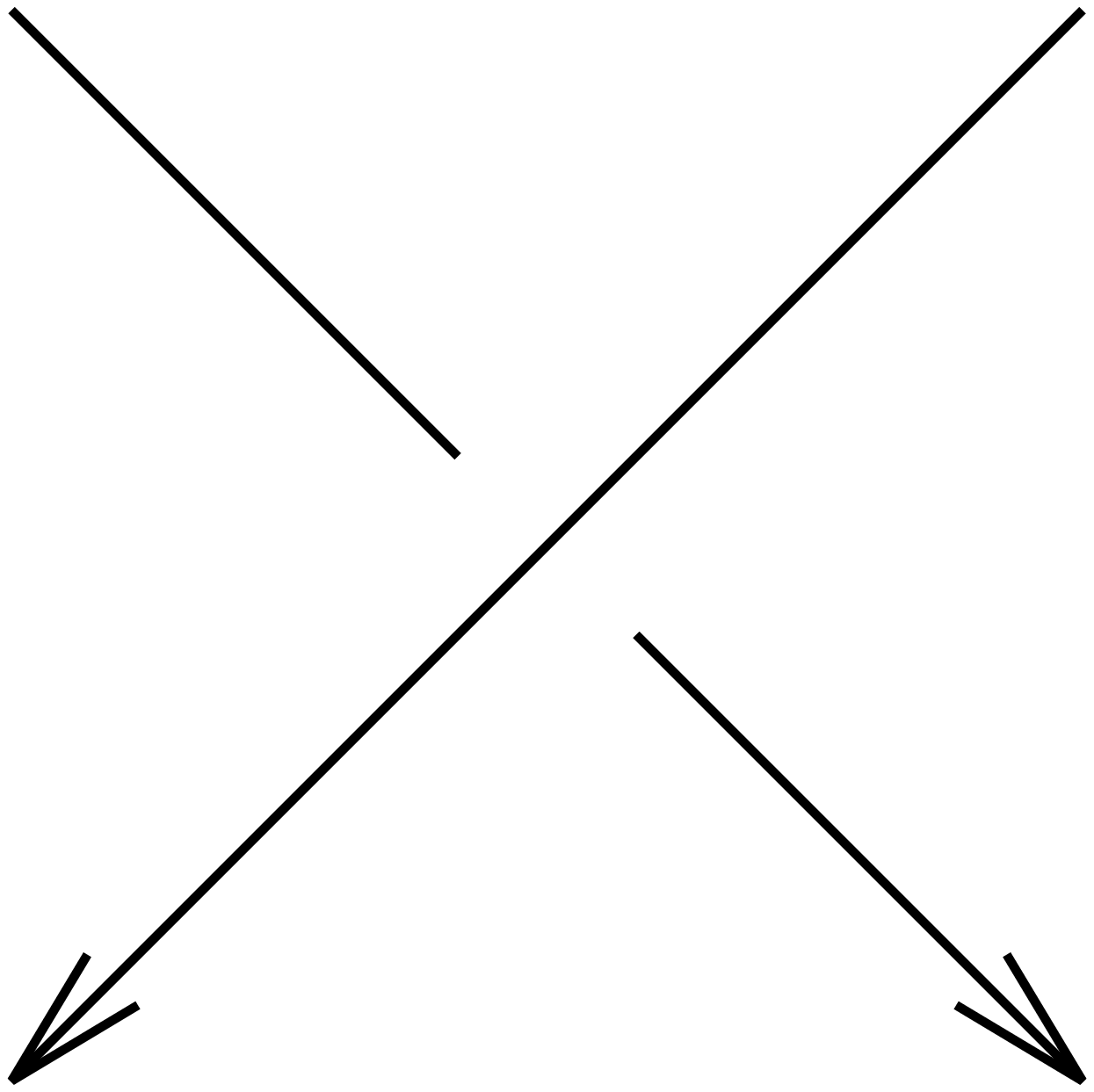}}
\medskip

\p for every triple of singular knots which differ at one 
crossing as indicated.
In particular, ${\cal K}_n$ can be viewed as a subspace of $\cal K$
for every $n$, by identifying any singular knot in ${\cal K}_n$
with the alternating sum of the $2^n$ knots obtained by
resolving its double points.
Hence, we have a subspace filtration
$$\ldots \subset {\cal K}_n\ldots
\subset{\cal K}_2\subset {\cal K}_1 \subset {\cal K}$$

\medskip
\p {\bf Definition 1.1.} {\sl A  Vassiliev  knot invariant of order $\leq n$
is a linear functional on the space ${\cal K}/ {{\cal K}_{n+1}}$.
The invariants of order $\leq n$ form a subspace ${\cal V}_n$
of ${\cal K}^*$, the annihilator of the subspace
${\cal K}_{n+1}\subset {\cal K}$. We will say that
an invariant $v$ is of {\it order n} if  $v$ lies in ${\cal V}_n$
but not in  
${\cal V}_{n-1}$. }
\medskip

Clearly, we have a filtration
$$ {\cal V}_0 \subset{\cal V}_1\subset {\cal V}_2 \subset \ldots$$

To continue we need to introduce some notation.
Let $D=D(K)$ be a diagram of
a knot $K$, and let ${\cal C}={\cal C}(D)=
\{ C_1,\ldots, C_m\}$ be a collection of disjoint non-empty sets of crossings
of $D$. Let us denote by $2^{\cal C}$
the set of all subsets of $\cal C$. Finally,
for an element $C\in 2^{\cal C}$ we will denote
by $D_C$ the knot diagram obtained from $D$
by switching the crossings in all sets contained in $C$. So, all together, we
can get $2^m$ different knot diagrams from the pair $(D,{\cal C})$. Notice that
each $C_i\in\cal C$ may contain more than one crossings.

\smallskip
\p {\bf Definition 1.2.} ([Gu])  {\sl Two knots  $K_1$ and $K_2$ are called 
{\it n-equivalent},
if $K_1$ has a knot diagram $D$ with the following property:
There exists ${\cal C}=\{C_1,\ldots, C_{n+1}\}$, a collection of $n+1$ disjoint
non-empty sets of crossings of $D$, such that $D_C$
is a diagram of $K_2$
for every non empty $C\in 2^{\cal C}$.
A knot $K$ which is {\it n-equivalent}
to the trivial knot will be called {\it n-trivial}.}
\medskip

By [Ya] every knot is {\it 1-trivial}.
The following theorem will be useful to us in the subsequent sections.

\medskip
\p {\bf Theorem 1.3.} ([Gu], [N-S]) {\sl
Two knots $K_1$ and $K_2$
are {\it n-equivalent} if and only if all of their Vassiliev invariants of 
order $\leq n$ are equal.
In particular, a knot $K$ is $n$-trivial
if and only if all its Vassiliev invariants of order
$\leq n$ vanish.}
\bigskip
\medskip

\centerline {2. {\B {\uppercase { Regular Seifert surfaces of a knot}}}}
\medskip
\p {\bf 2a. Generalities}
\medskip

Let $K$ be an oriented knot in $S^3$. A {\it Seifert surface} of $K$
is an oriented, compact, connected, bicollared surface
$S$, embedded in $S^3$ such that $\partial S=K$.

A {\it spine} of $S$ is a bouquet of circles $\Sigma \subset S$,
which is a deformation retract of $S$.

\medskip
\p {\bf Definition 2.1.} {\sl A Seifert surface $S$ of a knot $K$ is called
{\it regular} if it has a spine $\Sigma$ whose embedding in $S^3$,
induced by the embedding $S\subset S^3$, is
isotopic to the standard embedding of a bouquet of circles.
We will say that $\Sigma$ is a {\it regular spine} of $S$.}
\medskip

Let $\Sigma_n\subset S^3$, be a bouquet of $n$ circles
based at a point $p$. A regular projection of $\Sigma_n$ is a projection
of $\Sigma_n$ onto a plane with only transverse double points as
possible singularities. Starting from a regular projection of
$\Sigma_n$,
we can
construct an embedded compact oriented surface as follows: On the projection
plane, let
$D^2$ be a disc neighborhood of the basepoint $p$, which contains no
singular points of the projection. Then, $D^2$ intersects the
projection of $\Sigma_n$ in a bouquet
of $2n$ arcs and there are $n$ arcs outside $D^2$.
We first replace each of the arcs outside $D^2$ by a flat 
band with the original
arc as its core. Here a band being flat means we have an immersion when the
band is projected onto the plane. That is to say that the only singularities
the band projection has are these at the double points of the original arc
projection so that bands overlap themselves exactly when the arcs overcross
themselves. 

Let $S$ denote the surface obtained by the union of the disc
$D^2$ and these flat bands, to which some full twists are added if necessary. 
We say that $S$ is
a surface {\it associated} to the given regular plane projection of 
$\Sigma_n$. We 
will also
say that the surface $S$ is in a {\it disc-band form}. A
{\it band crossing} of $S$ is obviously defined, and they are in one-one 
correspondence
with crossings on the regular plane projection of $\Sigma_n$.
 We certainly have 
the
freedom to move the full twists added to the bands anywhere. 
So we assume that
all the twists of the band are moved near the ends of the bands.
 We may sometimes
abuse the notation by not distinguishing 
a band and its core and only take care of
the twists at the end of an argument.

Now let $S$ be a regular surface of genus $g$.
Pick a basepoint $p\in S$, and let
$\Sigma_n$, $n=2g$, be a 
regular spine of $S$ such that $p$ is the point on $\Sigma_n$ where all circles
in $\Sigma_n$ meet. Let $\gamma_1, \beta_1,\ldots,\gamma_g, \beta_g$
be the circles in $\Sigma_n$ oriented so that they form
a symplectic basis
of $H_1(S)$. Assume further that a disc neighborhood 
of $p$ in $S$ is chosen so that its intersection with $\Sigma_n$
consists of $2n$ arcs.

\medskip
\p {\bf Lemma 2.2.} 
{\sl Let $S$ be a regular Seifert surface,
with $\Sigma_n$ a regular spine. The embedding $\Sigma_n \subset S^3$
has a regular plane projection as shown in Figure 1 below,
where $b$ is a braid of index $2n$, such that
the regular Seifert surface $S$ is isotopic
to a surface associated to that projection of $\Sigma_n$. }
\medskip
\p {\it Proof.}
Let $W_n$ be a bouquet of $n$ circles, all based at a common point $q$. Then,
$\Sigma_n$ induces an embedding of $W_n$ in $S^3$.
Let us begin with a regular plane projection
of $W_n$, such that in a neighborhood
$D$ of $q$ in the projection plane,
the $2n$ arcs in $D\cap W_n$ are ordered and oriented
in the same way as the arcs of $\Sigma_n$
in the chosen disc neighborhood of $p$ in $S$.
Then, after a possible adjustment by adding some small kinks,
$S$ is isotopic to the surface associated to this
projection. Since $\Sigma_n$
is isotopic to the standard embedding of $W_n$ in $S^3$,
we may switch the arcs in $D$, so that the arcs of $\Sigma_n$
outside $D$ are isotopic to the standard embedding. We may then record these
switches by the braid $b$. \qed
\bigskip
\centerline{\epsfxsize=1.5in\epsfbox{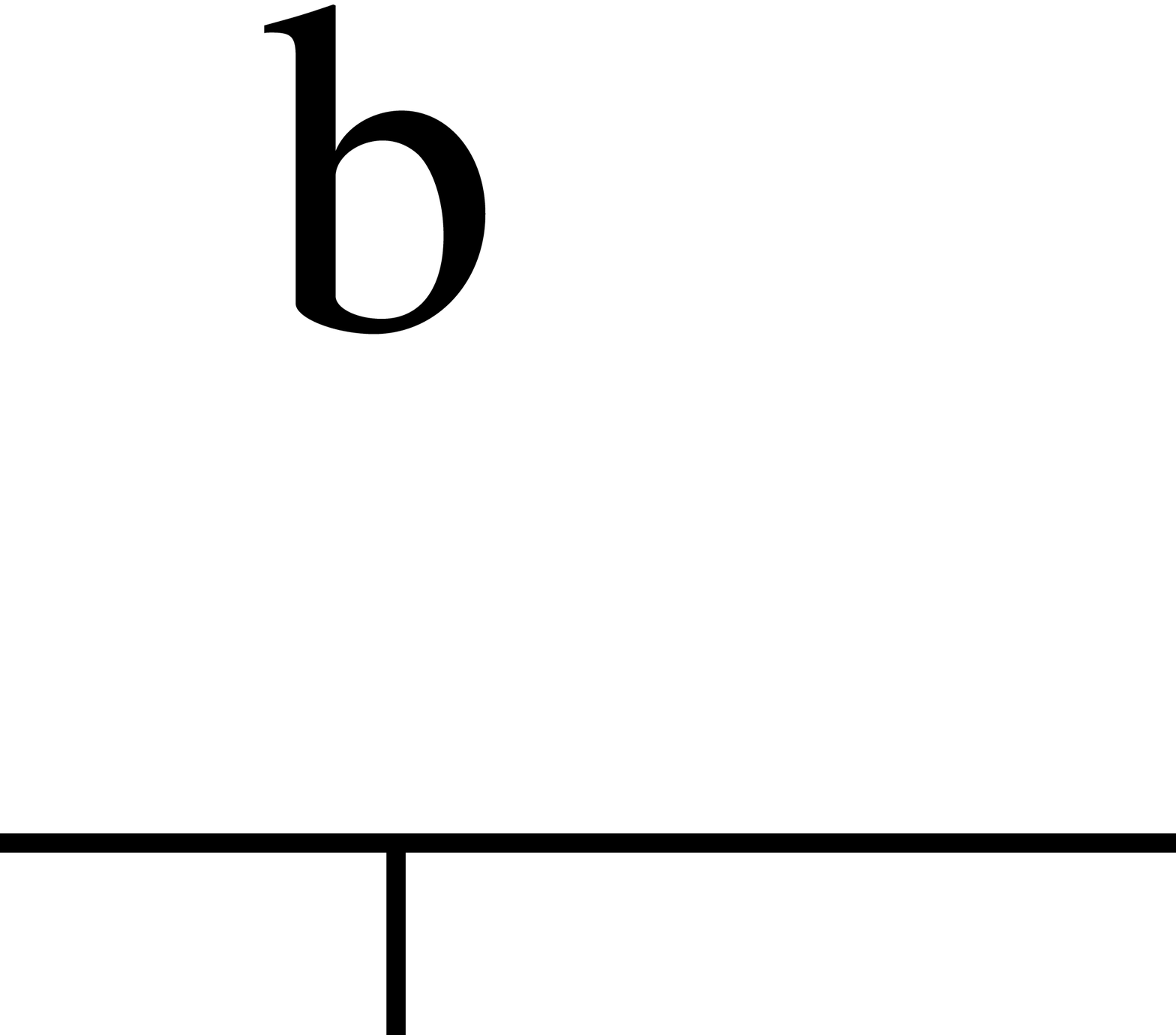} }
\medskip
\centerline {{\bf Figure 1.} {\msmall A projection of a regular spine
of a Seifert surface}}
\bigskip

\p {\bf 2b. Good position of a band in a regular spine}
\medskip

Let us consider $\R^3 \subset S^3$ and a decomposition
$\R^3=\R \times \R^2$, and take the factor $\R^2$ as a fixed
 projection plane $P$ 
from now on. Also, we will fix a coordinate
decomposition $(t,s)$ of $P$. Now let $l$ denote the $t$-axis on $P$,
and let $H_+=\{(t,s) \in P\, |\, s > 0\} $ 
and  $H_-=\{(t,s) \in P\, |\, s < 0\}$. 

To continue assume that $S$ is a regular Seifert surface and fix a
projection, $p: S \longrightarrow P$
in the disc-band form. Assume that the bands (their cores) of 
$S$ are all transverse to
$l$. Let
$B$ be a band of
$S$, and let
$\gamma$ be the core of
$B$. By a {\it sub-band} $B'$ of $B$, we mean a band on $B$ whose core $\gamma'$
is  a sub-arc of
$\gamma$.

\medskip
\p {\bf Definition 2.3.} {\sl We will say that a band $B$
is in {\it good position} with respect to the projection iff
the following conditions hold:

\p a) The band $B$ is flat.

\p b) For every band $A\neq B$,
the intersection $H_+ \cap A$ consists of a single sub-band with no
self-crossings, and these sub-bands are all disjoint.

\p c) All the self-crossings of $B$ occur in $H_+$,
and are all undercrossings (resp.
overcrossings).
Moreover, the intersection
$H_+ \cap B$ consists of finitely
many sub-bands $B_0,B_1,\ldots,B_k$ such that 

i) they have no self-crossings;

ii) the $B_j$'s, $j\neq0$, are disjoint with each other, and each crosses exactly once

under (resp. over) $B_0$, or one of the sub-bands in b). 

\p d) The crossings between $B$ and any other band that
occur in $H_-$ are all overcrossings (resp. undercrossings).}
\medskip

An example of a projection as described in Definition 2.3 
is shown in Figure 6, at the end of this section.

To continue, let $S$ be a regular Seifert surface and fix a
projection as described in Lemma 2.2.
Let $g$ be the genus and let $A_1, B_1, \ldots, A_g, B_g $
denote the bands of $S$. Moreover, let
$\gamma_1,\ \beta_1, \ \ldots, \ \gamma_g,\ \beta_g$
denote the cores of $A_1, B_1, \ldots, A_g, B_g $, respectively.
We orient the core curves so that they give a symplectic basis of $H_1(S)$.
Finally, let $x_1,y_1,\ldots,x_g,y_g$ be small linking 
circles of the bands such that

i) lk$(x_i, \gamma_j)=$ lk$(y_i, \beta_j)=\delta_{ij}$;

ii) lk$(y_i, \gamma_j)=$ lk$(x_i, \beta_j)=0$ and

iii) their projections on the plane $P$ are simple 
curves disjoint from each other. 

\p Clearly, $x_1,y_1,\ldots,x_g,y_g$ represent 
free generators of $\pi_1(S^3\setminus S)$.

\medskip
\p {\bf Lemma 2.4.} {\sl For every band $B$ of $S$,
there exists a projection of $S$
with respect to which $B$ is in good
position.}
\medskip

\p {\it Proof.} Let us start with the projection fixed
before the statement of the Lemma,
and let $l$ and $H_+$, $H_-$
be as before Definition 2.3.

Let
$\alpha_1,\hat\alpha_1,\ldots,\alpha_{2g-1},\hat \alpha_{2g}$
denote the hooks in Figure 1 on the top of the braid $b$. 
They are sub-bands of $A_1, B_1,
\ldots, A_g, B_g $, respectively. We move the projection of $S$ so that $l$
intersects
each
 of the hooks at exactly two points and we have that
the intersection $H_+ \cap p(S)$ is equal to
$\alpha_1\cup {\hat \alpha_1}\cup \ldots\cup 
\alpha_{2g-1}\cup{\hat \alpha_{2g}}$.
Thus the entire braid $b$ is left below $l$, in $H_-$.
To continue, we choose another horizontal line $l_0$ below $l$,
so that $b$ lies between
$l_0$ and $l$, and only the disc
part of the surface is left below $l_0$.
Finally, we draw more horizontal lines $l_1,l_2,\ldots,l_m=l$ 
such that the braid $b$
has exactly one crossing between $l_{i-1}$ and $l_i$ for $i=1,2,\ldots,m$.

Without loss of generality we may assume that $B=A_1$.

Observe that since $b$ is a braid, each 
band crossing of $A_1$ under some band $A$
can be slided all the way up,
by using the {\it finger moves} of Figure 2.
That is, we can slide a short sub-band of $A_1$, 
which is underneath $A$ at the crossing, up
following
$A$ until it becomes a small hook above $l$ undercrossing the hook of $A$. 

To isotope the band $A_1$ into good position, 
we start with the lowest undercrossing of
$A_1$ under, say some band $A$, between $l_{i-1}$ and $l_i$, for some $i$. 
We slide it up
above
$l$ and still call the resulting band $A_1$.
Now between $l_i$ and $l_{i+1}$, if there
is an undercrossing of $A$ in the original picture, we will have two 
new undercrossings
of the modified $A_1$. We slide these 
two new undercrossings of $A_1$ up, above
$l$, along the same way as we slide the undercrossing of $A$ between $l_i$ and
$l_{i+1}$ up above $l$. 
\bigskip
\centerline{\epsfxsize=3in\epsfbox{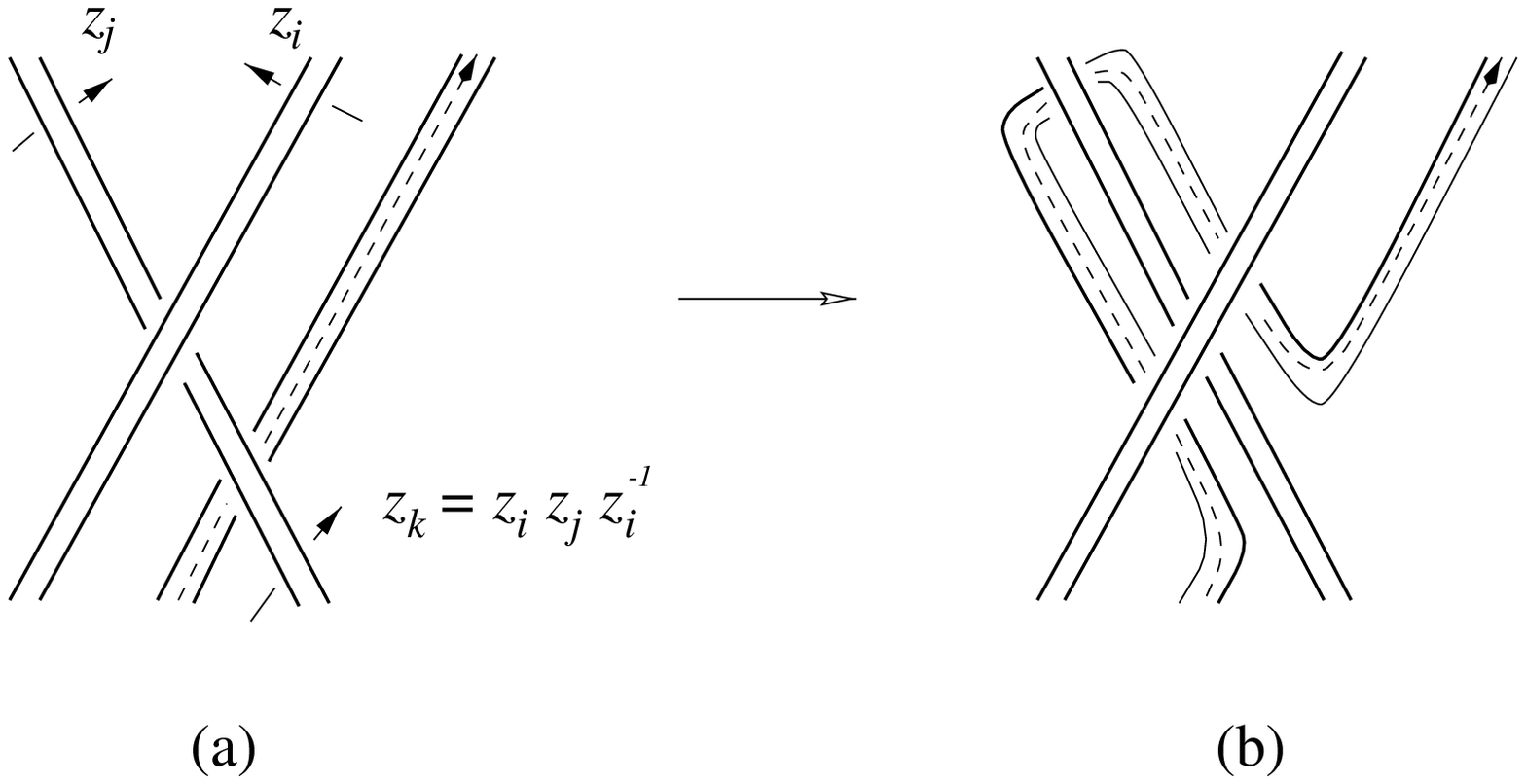} }
\medskip
\centerline {{ \bf Figure 2.} {\msmall Sliding an undercrossing across a band}}
\bigskip

Since $b$ has only finitely many crossings, this procedure will
slide all undercrossings of $A_1$ up above $l$, to make $A_1$
in good position.

\bigskip
\centerline{\epsfxsize=2.5in\epsfbox{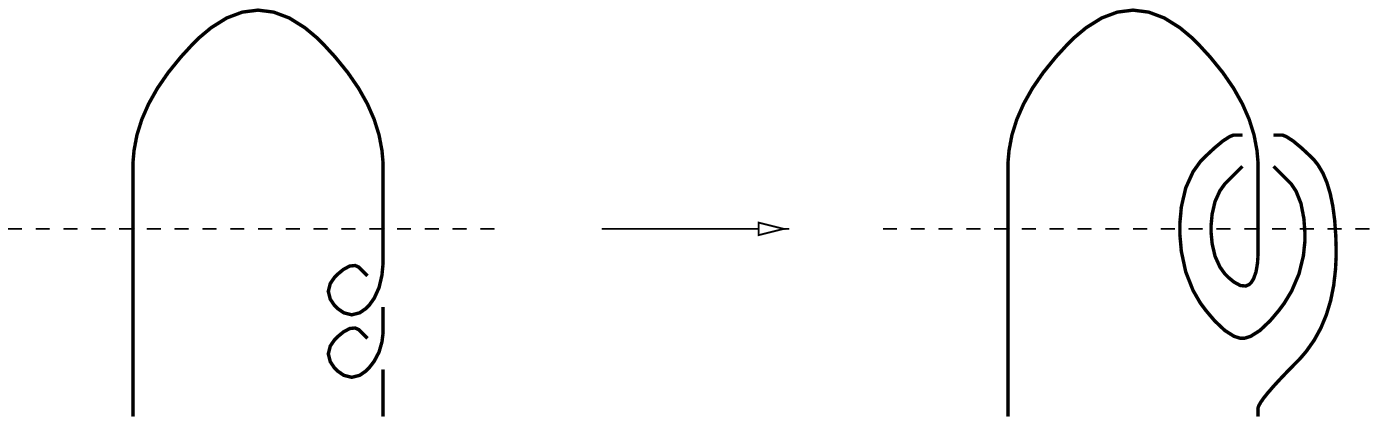} }
\medskip
\centerline {{\bf Figure 3.} {\msmall Twists on a band realized as kinks
and nested kinks}}
\bigskip

The
condition {\sl a)} of Definition 2.3 can also be satisfied by 
further isotopy which first moves the
twists on the band $B$ to a place around the line $l$ 
and then changes them to a
family of
\lq\lq nested kinks". See Figure  3.
\qed
    
\medskip
\p {\bf Remark 2.5.} Notice that a similar procedure can
be carried out for $B$ with overcrossings replacing undercrossings 
and vice versa.
\medskip

To continue let $K$ be a knot,
and let $S$ be a genus $g$ Seifert surface of $K$.
Let $N\cong S\times (-1,1)$  be a bicollar of $S$ in $S^3$,
such that $S\cong S\times \{0\}$, and let
$N^+= S\times (0,1)$ (resp. $N^-= S\times (-1,0)$).
For a simple closed curve $\gamma \subset S$, denote
$\gamma ^+=\gamma\times\{1/2\}\subset N^+$ and 
$\gamma^-=\gamma\times\{-1/2\}\subset N^-$.

Now let the projection
$p: S \longrightarrow P$, on the 
$(t,s)$-plane be as in Lemma 2.2. 
We denote by $z_1,z_2,\ldots,z_s$
the generators of $\pi_1(S^3\setminus S)$
arising from the Wirtinger presentation associated to the fixed
projection of $S$
(see for example [Ro]). The generators $z_i$ are in one to one correspondence
with the arcs of the projection
between two consecutive undercrossings. Moreover, every $z_i$
is a conjugate of one of the free generators
fixed earlier.

In Figure 2 we have indicated the 
Wirtinger generators
by small arrows under the bands.
The directions of the arrows are determined
by the fixed orientations
of the free generators
$x_1,y_1,\ldots,x_g,y_g$.
Notice that these free generators can be chosen as a part of the
Wirtinger 
generators
corresponding to the hooks $\alpha_1, {\hat \alpha_1} , \ldots,
\alpha_{2g-1},
 {\hat
\alpha_{2g}}$, respectively. 

\medskip
\p {\bf Lemma 2.6.} (The geometric rewriting) {\sl Let $\gamma \subset S$ 
be a simple closed curve represented by
the core of a band $B$ of $S$, and let
$\gamma^{\epsilon}$  ($\epsilon=\pm$) be one of the push off's of
$\gamma$. Moreover, let $W=z_{i_1}\cdots z_{i_m}$ be
a word representing $\gamma^{\epsilon}$
in terms of the Wirtinger generators
of the projection, and let $W'=W'(x_i, y_j)$ be the word obtained from $W$
by expressing each $z_{i_r}$ in terms of the free generators.
Then, there exists a projection
$$p': S \longrightarrow P$$
\p which is obtained from
$p$ by isotopy, and such that
every letter in $W'$ is realized by an undercrossing of $B$ with one of the
hooks above $l$.}
\medskip

\p {\it Proof.} Suppose that 
${\epsilon}=+$. Assume first that $B$ is flat. Then every $z_{i_r}$
in $W$ can be realized by an undercrossing of $B$ with another
band or itself.
We claim that the projection obtained in Lemma 2.4
has the desired properties. 

To see that, let us begin with
three Wirtinger generators $z_i,zj,z_k$,
around a crossing of the projection as show
in Figure 2 (a). 
Then we have $z_k=z_iz_jz_i^{-1}$. Notice that $\gamma^+$,
which is drawn by the dotted arrow,
picks up $z_k$ at the undercrossing in (a).
After performing a {\it finger move}
in (b), each of $z_i$, $z_j$, $z_i^{-1}$ is realized 
by an undercrossing.
Thus the geometric equivalent of replacing each $z_{i_j}$,
in $W=z_{i_1} \cdots z_{i_m}$, by its expression
in terms of the free generators
is to slide
an undercrossing of
$B$ until it reaches the appropriate hook.
Now the desired conclusion follows easily from this observation.

If $B$ is not flat, we first assume that the twists are all near 
one of the ends of $B$. Let the Wirtinger generator near that end of $B$ be 
$z_0$. Then $W=z_0^k
z_{i_1}\cdots z_{i_m}$, and every $z_{i_r}$ can be realized by an 
undercrossing of $B$ with another band or itself. The previous argument 
still works in this case for
each $z_{i_j}$. For $z_0^k$, we may first move the twists up to a place 
around $l$ and then replace the twists by a family of nested kinks, like in 
the last part of the proof of Lemma 2.4. This finishes the proof of Lemma
2.6.
\qed

\medskip
\p {\bf 2c. Lower central series and curves on surfaces}
\medskip

For a group $G$ let $[G, G]$ denote the commutator subgroup
of $G$. The lower central series
$\{ G^{(m)}\}_{m\in \N}$, of $G$ is defined
by
$G^{(1)}=G$ and
$$G^{(m+1)}=[G^{(m)}, G]$$
\p for $m\geq 1$.
\medskip
We begin by recalling some commutator identities that
will be useful to us later on. See [KMS].

\medskip
\p {\bf Proposition 2.7.} (Witt-Hall identities) 
{\sl Let $G$ be a group and let $k$,
$m$ and $l$ be positive integers. Suppose
that $x\in G^{(k)}$, $y\in G^{(m)}$
and $z\in G^{(l)}$. Then

\p a)\ \ \ $[G^{(k)},\ G^{(m)}] \subset G^{(k+m)}$
or $xy \equiv yx$ mod $G^{(k+m)}$

\p b)\ \ \ $[x, zy] = [x, z]\  [x, y] \ [ [y, x], z] $

\p c)\ \ \ $[xy, z] = [y, z]\ [ [z, y], x]\ [x, z]  $

\p d)\ \ \ $[x, [y, z]]\  [y, [z, x]]\  [z, [x, y]]\equiv 1$
 mod $G^{(k+l+m+1)}$

\p e)\ \ \  If $g \equiv g'$ mod $G^{(k)}$
and $y\in G^{(m)}$ then
$[g, y] \equiv [g', y]$ mod $G^{(k+m)}$
and
$[y, g] \equiv [y , g']$ mod $G^{(k+m)}$.}
\medskip

Let $F$ be a free group of finite rank
and let ${\cal A}=\{ a_1, \ldots, a_k\}$
be a set of (not necessarily free) generators of $F$.
Let $a$ be an element in $F$
and let $W=W_a({\cal A})$ be a word
in $a_1, \ldots, a_k$
representing $a$. Think of $W$ as given as a list of spots in which we
may deposit
letters $a_1^{\pm1},a_2^{\pm1},\ldots,a_k^{\pm1}$. Now let
${\cal C}={\cal C}(W)=
\{ C_1,\ldots, C_m\}$ be a collection of disjoint non-empty sets of spots (or letters)
in $W$. Let us denote by $2^{\cal C}$
the set of all subsets of $\cal C$. Finally,
for an element $C\in 2^{\cal C}$ we will denote
by $W_C$ the word obtained from $W$
by substituting the letters in all sets contained in $C$ by $1$.

\medskip
\p {\bf Definition 2.8.} ([N-S]) {\sl The element
$a\in F$ is called {\it n-trivial},
with respect to $\cal A$,
if it has a word presentation $W=W_a({\cal A})$
with the following property:
There exist a collection of $n+1$ disjoint non-empty sets of letters, 
say ${\cal C}=\{C_1,\ldots, C_{n+1}\}$, in $W$ such that $W_C$
represents the trivial element
for every non-empty $C\in 2^{\cal C}$.}
\medskip

We will say that $a\in F$ is {\it n-trivial}
if it is {\it n-trivial} with respect to a set of generators.
The following lemma shows that this definition depends neither on the word
presentation nor the set of generators used.

\medskip
\p {\bf Lemma 2.9.} {\sl If $a\in F$
is {\it n-trivial}
with respect to a generating set $\cal A$,
then it is {\it n-trivial}
with respect to every generating set of $F$.}
\medskip

\p {\it Proof.} Let
$W=W_a({\cal A})=a_{i_1}a_{i_2} \ldots a_{i_s}$ be a word for $a$
satisfying the properties in Definition 2.8
and let ${\cal A}'$ be another set of generators for $F$.
By expressing each $a_{i_j}$ as a word of elements
in ${\cal A}'$ we obtain a word $W'=W'_a({\cal A}')$
which satisfies the requirements of n-triviality
with respect to ${\cal A}'$. \qed

\medskip
\p {\bf Lemma 2.10.} {\sl If $a$ lies in $F^{(n+1)}$,
then it is {\it n-trivial}.}
\medskip

\p {\it Proof.}  Observe that a basic commutator
$[a,\ b]=aba^{-1}b^{-1}$ is 1-trivial by using
${\cal C}= \{ \{a, a^{-1}\}, \ \{b, b^{-1}\}\}$,
and induct on $n$. \qed

Clearly, we do not change the $n$-triviality of a word
by inserting  a {\it canceling
pair} $xx^{-1}$ or $x^{-1}x$, where $x$ is a generator. 
We will use the following definition
to simplify the exposition.

\medskip
\p{\bf Definition 2.11.} {\sl A {\it simple commutator} is a word in the form
of $[A,x^{\pm1}]$ or $[x^{\pm1},A]$ where $x$ is a generator and $A$ is a
simple commutator of shorter length. A {\it simple quasi-commutator} is a word
obtained from a simple commutator by finitely many insertions of
canceling pairs.}
\medskip

By Proposition 2.7, any word representing an element 
in $F^{(n)}$ can be changed
to a product of simple quasi-commutators of length $\geq n$ by finitely many
insertions of canceling pairs. A simple quasi-commutator of length
$>n$ is clearly $n$-trivial.

To continue, let $S$ be a regular Seifert surface
of a knot $K$. For a loop 
$\alpha\subset S^3 \setminus S$, 
we will denote by $[\alpha]$ its homotopy class in
$\pi=\pi_1( S^3 \setminus S)$.

Suppose that $S$, $\gamma$, $B$, $p': S \longrightarrow P$ and
$[\gamma^{\epsilon}]=W'(x_i,y_i)=W'$ are as in the statement of Lemma 2.6.
Suppose $\delta$ is a sub-band of $\gamma$ and $[\delta^{\epsilon}]=W''$
is a sub-word of $W'$. Assume that $W''$ represents an element in $\pi^{(n)}$.

\medskip
\p {\bf Lemma 2.12.}
(The geometric realization)
{\sl There exists   
a projection $p_1:S\rightarrow P$ with the following properties:

\p i) $p_1(S)$ is obtained from $p'(S)$
by a finite sequence of band Reidermeister moves of type II; 

\p ii) $B$ is in good position with respect to the new projection; 

\p iii) the word $W^*=W^*(x_i,y_i)$ one reads out from $\delta$ 
(with respect to the new projection 
$p'$), 
by picking up one letter for each
crossing of $\delta$ underneath the hooks, 
is a product of simple quasi-commutators of
length $\geq n$.}
\medskip

\p {\it Proof.} Since $W''$ represents 
an element in $\pi^{(n)}$,
we may change $W''$ to a product of simple
quasi-commutators by finitely many
insertions of canceling pairs. Such an insertion of a cancelling pair can 
be realized geometrically by a finger move (type II Reidermeister move). 
We will create a region in the projection plane to perform such a finger move. 
This region is a horizontal long strip below the line $l$ and its 
intersection with $p'(S)$ consists of vertical straight flat bands. 
See Figure 4.

\bigskip
\centerline{\epsfxsize=3.5in\epsfbox{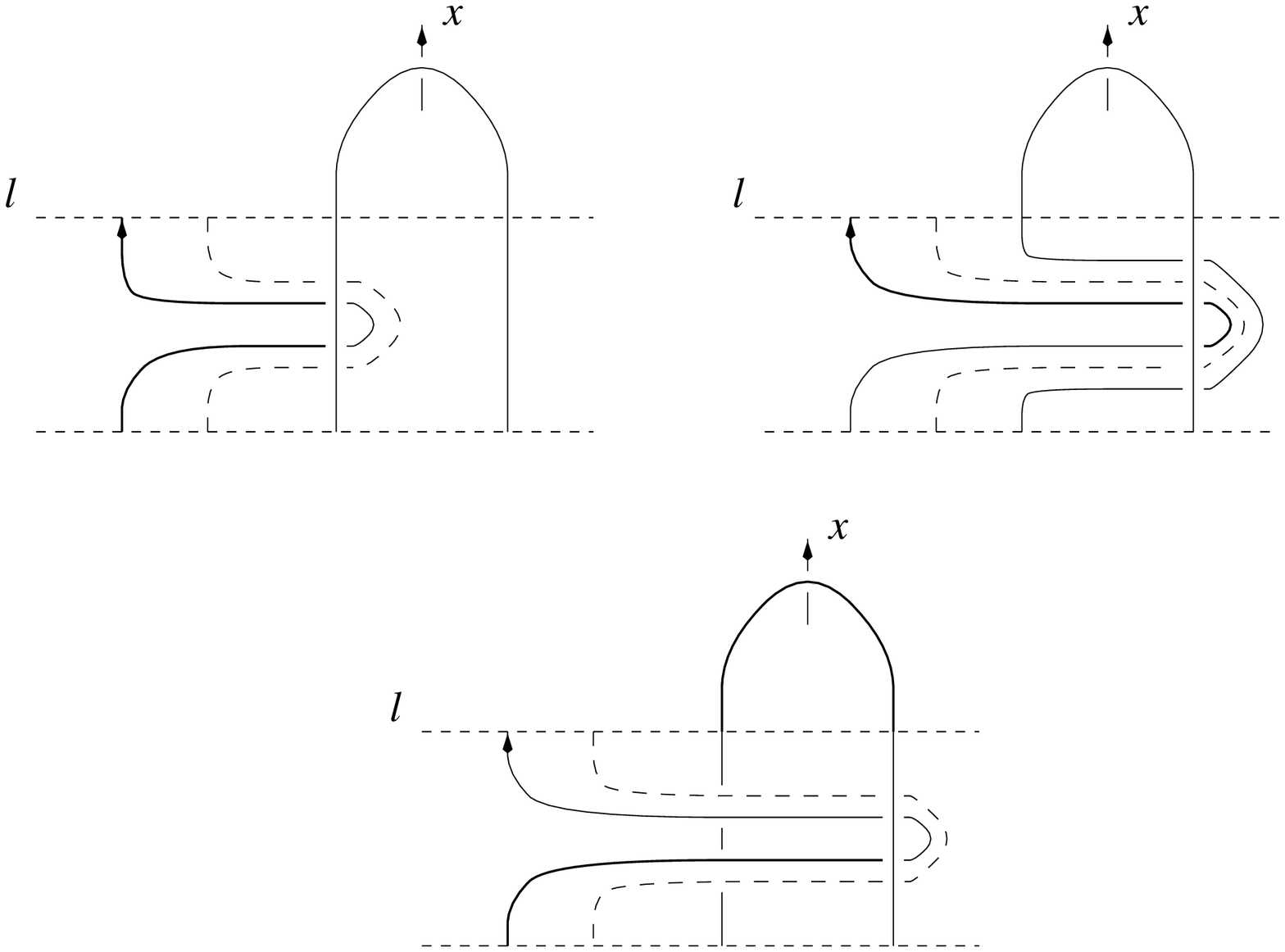} }
\medskip
\centerline {{ \bf Figure 4.}
{\msmall Realizing an insertion of $x^{-1}x$ or $xx^{-1}$ by isotopy}}
\bigskip
 As shown in Figure 4, there are two situations corresponding to
insertions of $x^{-1}x$ or  $xx^{-1}$. In one of the cases, we shall 
either let the finger go over one of the vertical flat bands connected to 
the $x$-hook (in the case that the $x$-hook does not belong to $B$) or
push that vertical flat band along with the finger move (in the case
that the $x$-hook belongs to $B$).
Furthermore, if some vertical flat bands belonging to $B$ block the
way of the finger move, we will make more insertions by pushing these
vertical flat bands along with the finger move. Finally, with all these done,
we may easily modify the projection further to make $B$ still in good position
and ready to do the next insertion. \qed
\medskip
\medskip
\p {\bf 2d. An example illustrating the good position of a band} 
\medskip
In Figure 5, we show an example of a 
regular Seifert surface of genus one.  
\medskip
\centerline{\epsfxsize=1in\epsfbox{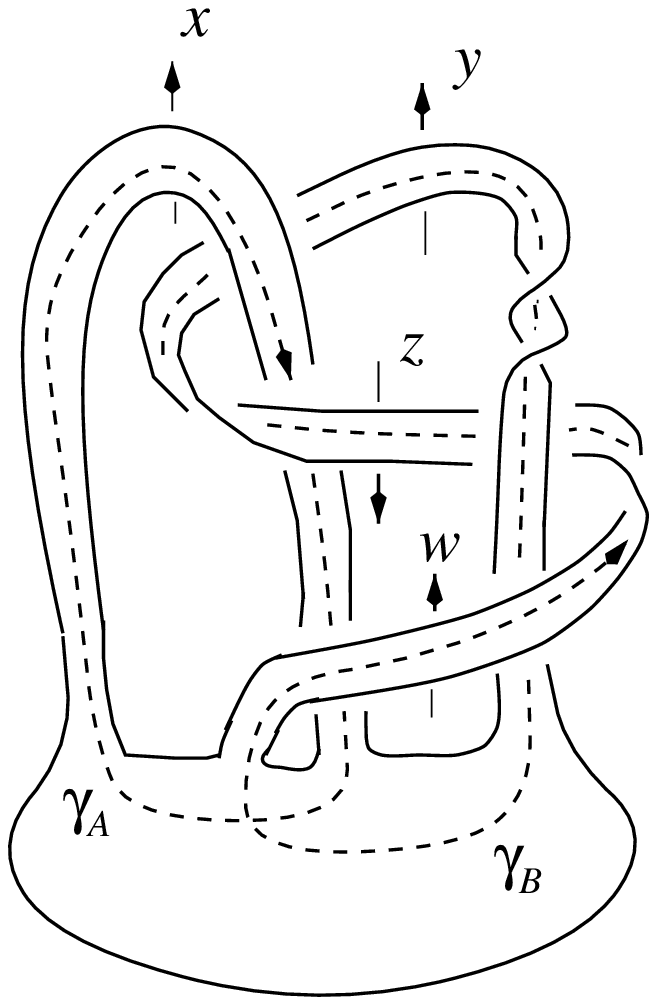} }
\medskip
\centerline {{ \bf Figure 5.}\ {\msmall A 2-hyperbolic surface }}
\medskip
The cores of the bands 
$A$ and $B$
of the surface $S$ have been
drawn by the dashed,
oriented curves $\gamma_A$ and $\gamma_B$,
respectively. The fundamental group
$\pi=\pi_1(S^3\setminus S)$ is freely generated by $x$ and $y$. We have
$[\gamma_A^+]\in\pi^{(3)}$, where $\gamma_A^+$ is the push-off of $\gamma_A$ 
along the positive normal vector of the surface $S$ 
pointing upwards the projection plane. In fact, from the 
Wirtinger presentation obtained from the given projection
we have $[\gamma_A^{+}]=[zw^{-1}]=[[x, y], y^{-1}]$. 
Such a surface  will be called {\it 2-hyperbolic} in 
Definition 3.1.  
Now we modify the projection of $S$, so that
$A$ is in good position. The resulting
projection is shown
in Figure 6.
Here we have only drawn
the cores of the bands. The solid (dashed, resp.)
arc corresponds to the 
band $A$ ($B$, resp.)
of Figure 6. 
The word we read out when traveling along the solid arc,
one letter for each crossing underneath the hooks, is exactly
$W=[[x, y],y^{-1}]$. 
The reader can verify that the two crossings marked by $``\%"$
on the left side of the $y$-hook and the two crossings marked
by $``\ast"$
can be used to trivialize the band $A$ in three
ways.
\bigskip
\centerline{\epsfxsize=2.5in\epsfbox{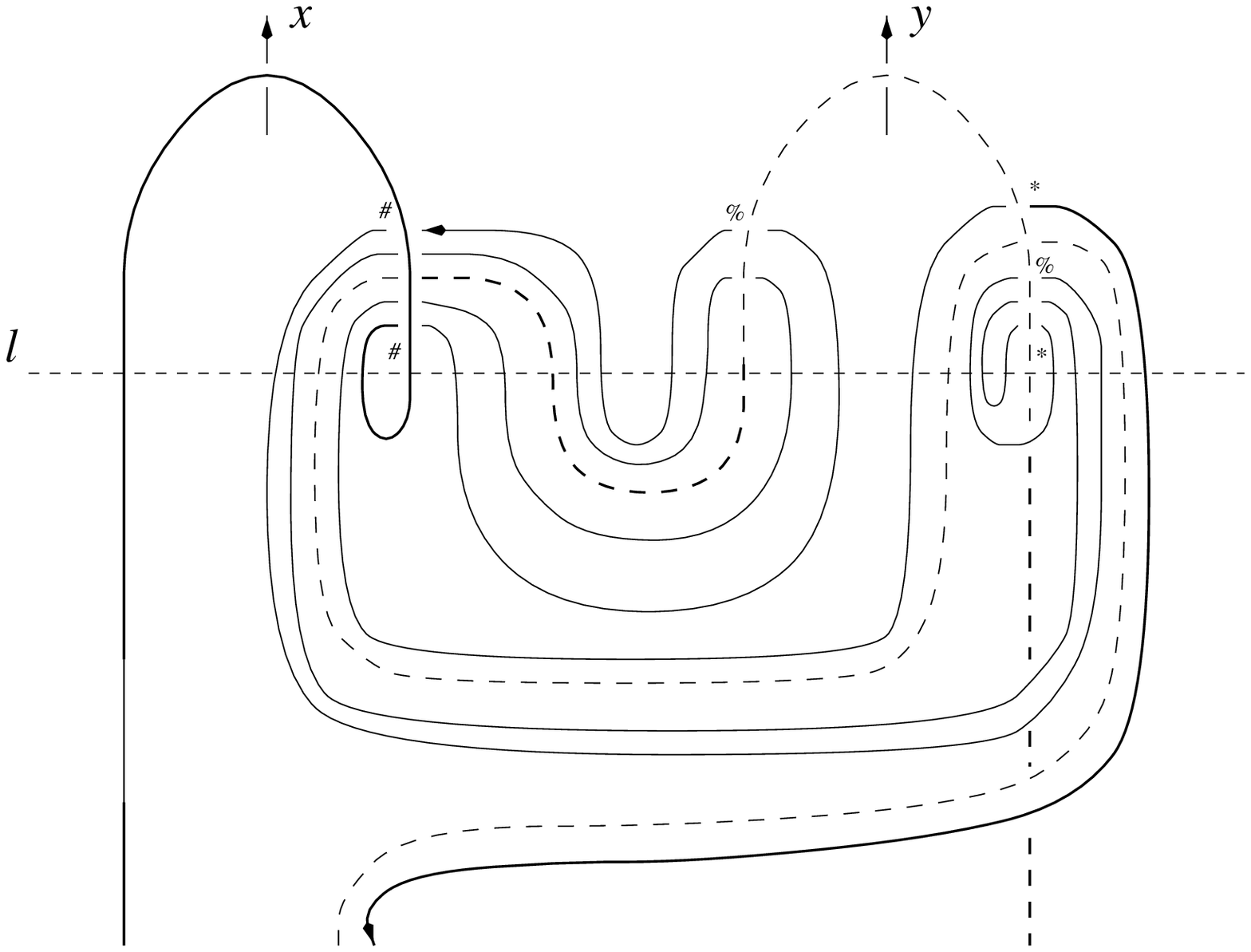} }
\medskip
\centerline {{ \bf Figure 6.}\ {\msmall The band 
$A$ in good position}}
\bigskip

\medskip
\p {\bf Remark 2.13.} There is an obvious collection $\cal C$
of three sets 
of letters of $W$, so that $W$ becomes a trivial word whenever
we delete letters
in a non-empty $C\in 2^{\cal C}$ from $W$ (see proof of Lemma 2.10).
The projection of Figure 6 has the property
that every letter in the word $W=[[x, y],y^{-1}]$ 
is realized by a band crossing. So we obtain a collection of sets of crossings,
also denoted by $\cal C$.
However, as the reader can verify,
the image of $\gamma_A$
on the surface $S_C$, obtained from $S$ by switching the crossings in $C$,
will not always be homotopically trivial in $S^3\setminus S_C$.
\bigskip
\medskip

\input ntmacros.tex
\centerline {{ 3. {\B \uppercase {commutator half bases and 
Vassiliev invariants}}}}
\medskip
In this section we undertake the study of
{\it  regular} Seifert surfaces, whose complement
looks, modulo the first 
${n+1}$ terms of the lower central series
of its fundamental group, like the complement
of a null-isotopy. Our main goal is to show that
the existence of such a surface for a knot $K$
forces its Vassiliev invariants
of certain orders to vanish.
\medskip
\p {\bf 3a. Definitions and the main result } 
\medskip
Before we are able to state our main result in this section
we need some notation and terminology.
Let $K$ be a knot in $S^3$ and let $S$ be a Seifert surface of $K$,
of genus $g$. Throughout this paper a {\it basis}
of $S$ will be a collection of
$2g$  non-separating simple 
closed curves $\{ \gamma_1, \beta_1, \ldots, \gamma_g, \beta_g \}$
that represent a symplectic basis of $H_1(S)$.
That is we have
${\cal I} (\gamma_i, \ \gamma_j)={\cal I}(\beta_i, \ \beta_j)=
{\cal I}(\beta_i, \ \gamma_j)=0$,
for $i\neq j$, and ${\cal I}(\gamma_i, \ \beta_i)=1$, where ${\cal I}$
denotes the intersection form on $S$.
Each of the collections $\{ \gamma_1, \ldots, \gamma_g \}$
and $\{ \beta_1, \ldots, \beta_g \}$ will
be called {\it a half basis}. 

To continue let $\pi={\pi_1}(S^3\setminus S)$.
For a basis  ${\cal B}= \{ \gamma_1, \beta_1, \ldots, \gamma_g, \beta_g \}$
of $H_1(S)$
let ${\cal B}^{*}= \{ x_1, y_1,
\ldots, x_g, y_g \}$ denote elements in $\pi$ representing the dual basis 
of $H_1(S^3 \setminus S)$.

For a subset $\cal A $ of $\cal B$, let  
${\cal G}_{\cal A}$ denote  the normal subgroup of $\pi$
generated by the subset of ${\cal B}^{*}$
corresponding to $\cal A$. Moreover, we will denote by
$\pi_{ \cal A}$ (resp.  $\phi_{\cal A}$)
the quotient $\pi / {\cal G}_{\cal A}$
(resp. the quotient  homomorphism
$\pi \longrightarrow \pi / {\cal G}_{\cal A}$).
Finally, $\pi_{ \cal A}^{(m)}$
will denote the $m$-th term of the lower
central series of $\pi_{ \cal A}$. For the following definition
it is convenient to allow $\cal A$ to be the empty set
and have $\pi_{ \cal A}=\pi$.
\medskip
\p {\bf Definition 3.1.} {\sl Let $n\in \N$.
A regular Seifert surface $S$ is called {\it n-hyperbolic},
if it has a half basis $\cal A$ 
represented by circles in a regular spine $\Sigma$
with the following property:
There is an ordering, $\gamma_1, \ldots, \gamma_g$,
of the elements in $\cal A$ such that either
$\phi_{{\cal A}_{i-1}}([
\gamma_i^+])$ or $\phi_{{\cal A}_{i-1}}([\gamma_i^-])$ lies in $\pi_{{\cal A}_{i-1}} ^{(n+1)}$.
Here  ${\cal A}_k=
 \{ x_1, y_1, \ldots, x_k, y_k
\}$
for $k
=1, \ldots, g$  and ${\cal A}_0$
is the empty set.
The boundary
of such a surface
will be called an {\it n-hyperbolic} knot.}
\medskip
\p {\it Remark.} A particularly interesting case of Definition 3.1 is when
the surface $S$ contains a half basis $\gamma_1, \ldots, \gamma_g$
such that $[\gamma_i^+]$ or $[\gamma_i^-]$ lie in $\pi^{(n+1)}$.
\medskip
In order to state our main result in this section
we need some notation.
For $m\in \N$, let $q(m)$ be the quotient of division
of $m$ by six (that is $m=6q(m)+r_1$, $0\leq r_1 \leq 5$).
Let the notation be as in  Definition 3.1.
For $i=1, \ldots, g$,
let $x_i$ denote the free generator of $\pi$ that is
dual to
$[\gamma_i]$. Let
$l_i$ denote the number of distinct elements
in
$\{ x_1, y_1, \ldots, x_g, y_g \}$, that are different than
$x_i$ and whose images under $\phi_{{\cal A}_{i-1}}$
appear in a (reduced) word, say $W_i$,
representing $\phi_{{\cal A}_i}([\gamma_i^+])$ or $\phi_{{\cal A}_i}(
[\gamma_i^-])$. Write $W_i$ as a product,
$W_i=W_i^1 \ldots W_i^{s_i}$, of elements
in $\pi_{{\cal A}_{i-1}} ^{(n+1)}$ and partition
the set $\{ W_i^1, \ldots W_i^{s_i}\}$
into disjoint sets, say ${\cal W}_i^1, \ldots {\cal W}_i^{t_i}$
such that: i) $k_i^1+\ldots + k_i^{t_i}=l_i$, where $k_i^j$
is the number of distinct elements in ${\cal A}_i$ involved
in ${\cal W}_i^j$ and ii) for $a\neq b$, the sets
of elements from ${\cal A}_i$ appearing
${\cal W}_i^a$ and ${\cal W}_i^b$ are disjoint.
Let $$k_i ={\rm min} \{k_i^1, \ldots, k_i^{t_i}\}.$$
\p and
let
$$q_{\gamma_i}= \cases {q(n+1), &if $n < 6k$\cr
\phantom{a} & \phantom{a} \cr
 k_i+ \left[{
{\displaystyle {\rm log}_2{
{{n+1-6k_i}\over 6}}}}
\right], &if $n \geq 6k$.\cr }
 $$                                          
Notice that
$$q(n+1)> {{n+1}\over 6} -1={{n-5}\over 6}>{\rm log}_2(
{{n-5}\over 6}).$$
Also, since $ab\geq a+b$ if $a,b>1$, we have
$$k_i+ {\rm log}_2({{n+1-6k_i}\over 6})
>{\rm log}_2k_i+ {\rm log}_2({{n+1-6k_i}\over 6})
>{\rm log}_2({{n+1}\over 36}).
$$
Thus, for $n>5 $, we have
$$q_{\gamma_i}>{\rm log}_2({{n-5}\over 72}).$$

We define $l(n)$ by
$$l(n, \ S)= {\rm min} \{ q_{\gamma_1}-1, \ldots, q_{\gamma_g}-1\},$$
\p and
$$l(n)={\rm min} \{ l(n, \ S)|\  S \ is\  n{\rm -hyperbolic}\}.$$
We can now state our main result in this section, which is:
\medskip
\p {\bf Theorem 3.2.} {\sl If $K$ {\it n-hyperbolic}, for some $n\in \N$,
then $K$
is at least {\it l(n)-trivial}. Thus, all the Vassiliev invariants
of $K$ of orders $\leq l(n)$ vanish.}
\medskip 
>From our analysis above, we  see that 
${\displaystyle {l(n) > {\rm log}_2({ n-5 \over 144}}})$
and in particular
${\displaystyle {\rm lim}_{n \to \infty }l(n) =\infty}$.
Thus, an immediate Corollary of Theorem 3.2 is:
\smallskip
\p {\bf Corollary 3.2*.} {\sl If $K$ {\it n-hyperbolic}, for all
 $n\in \N$,
then all its Vassiliev invariants vanish.}
\medskip
Assume that $S$ is
in disc-handle form as described in Lemma 2.2
and that the cores of the bands form a symplectic basis of $H_1(S)$.
Moreover, assume that
the curves $\gamma_1, \ldots, \gamma_g$,
of Definition 3.1 can be realized by half of these cores.
Let $\beta_1, \ldots, \beta_g$ denote the cores of the other half bands
and let $D=D(K)$ denote the knot diagram of $K$, induced by our projection
of the surface.
Also we may assume that the dual basis
$\{ x_1, y_1,
\ldots, x_g, y_g \}$ is represented by free generators
of $\pi$, as before the statement of Lemma 2.4.

To continue with our notation , let $C$ be a collection of band crossings
on the projection of $S$. We denote by $S_C$ (resp. $D_C$)
the Seifert surface (resp. knot diagram) obtained from $S$ (resp. $D$) 
by switching all crossings in $C$, simultaneously. For a simple curve 
$\gamma \subset S$ (or an arc 
$\delta \subset \gamma$), we will denote by $\gamma_C$ 
(or $\delta_C$) the image of 
$\gamma$ (or $\delta$)  on $S_C$. 

Let $\gamma$ be the core of a band $B$ in 
good position and suppose that it is decomposed into a union of
sub-arcs $\eta\cup\delta$ with disjoint interiors, such that 
the word, say $W$, represented by 
$\delta^+$ (or $\delta^-$) 
in $\pi=\pi_1(S^3\setminus S)$ lies in  $\pi^{(m+1)}$. 
Let $x$ be the generator of $\pi$ corresponding
to $B$ and let
$l$ denote the number of distinct
free generators, different than $x$,
appearing in $W$.
Write $W$ as a product,
$W=W_1 \ldots W_s$, of commutators
in $\pi^{(n+1)}$ and partition
the set $\{ W_1, \ldots W_s\}$
into disjoint sets, say ${\cal W}_1, \ldots {\cal W}_t$
such that: i) $k_1+\ldots + k_{t}=l$, where $k_j$
is the number of distinct generators involved
in ${\cal W}_j$ and ii) for $a\neq b$, the sets
of generators appearing
${\cal W}_a$ and ${\cal W}_b$ are disjoint.
Let $$k ={\rm min} \{k_1, \ldots, k_t\}.$$
We define $$q_{\delta}= 
\cases {q(n+1), &if $n < 6k$\cr
\phantom{a} & \phantom{a} \cr
k+ \left[ {\displaystyle {\rm log}_2({
{{n+1-6k}\over 6})}}\right], &if $n \geq 6k$.\cr }
$$
The proof
of Theorem 3.2 will be seen to follow from the following Proposition.
  
\medskip
\p {\bf Proposition 3.3.} {\sl  Let $\gamma$ be the core of a band $B$ in 
good position and suppose that it is decomposed into a union of
sub-arcs $\eta\cup\delta$, such that 
the word represented by 
$\delta^+$ (or $\delta^-$) 
in $\pi=\pi_1(S^3\setminus S)$ lies in  $\pi^{(m+1)}$. 
Suppose that the word represented
by $\eta^{+}$ (or  $\eta^{-}$)
is the identity.

Let $K'$ 
be the boundary of the surface obtained from $S$ by replacing 
the sub-band of $B$ corresponding to $\delta$ with a straight flat ribbon
segment $\delta^{*}$,
 connecting the endpoints of $\delta$ and above (resp. below) the remaining 
diagram. Then $K$ and $K'$ are at least $l_{\delta}$-equivalent,
where $l_{\delta}=q_{\delta}-1$.}
\medskip

The proof of Proposition 3.3 will be divided into several steps,
and occupies all of $\S 3$.
Without loss of generality we will work with
$\delta^{+}$ and $\gamma^{+}$. 
In the course of the proof we will see
 that we may choose
the collection of  sets of crossings
$ \cal C$, required in the definition of $l_{\delta}$-equivalence, 
to be band crossings
in a projection of $S$. Moreover,
for every non-empty $C\in 2^{\cal {C}}$,
$\delta_C$ will be shown to be  isotopic to 
a straight arc, say  $\delta^{*}$, as in the statement above.
Here $2^{\cal {C}}$ is 
the set of all subsets of $\cal C$.

In the rest of this paragraph, let us
assume Proposition 3.3 and show how Theorem 3.2 follows from it.
\medskip
\p {\it Proof of Theorem 3.2}. The proof
will be by induction on the genus
$g$ of the surface $S$. If $g=0$
then $K$ is the trivial knot and there is nothing to prove.
For $i=1,\ldots, g$, let $A_i$  denote the band of $S$
whose core corresponds to $\gamma_i$, and let $B_i$
be the dual band.

By Definition 3.1 we have a band $A_1$, such
that the core $\gamma$ satisfies the assumption of Proposition 3.3.
We may decompose
$\gamma$
into a union of
sub-arcs $\eta\cup\delta$ with disjoint interiors such that
the word represented by 
$\delta^+$ (resp.
 $\eta^{+}$)  in $\pi=\pi_1(S^3\setminus S)$ lies in  $\pi^{(n+1)}$
(resp. is  the empty word).
Let $K'$ 
be a knot obtained from $K$ by replacing 
the sub-band of $B$ corresponding to $\delta$ with a straight flat ribbon
segment $\delta^{*}$,
connecting the endpoints of $\delta$ and above the remaining 
diagram, and let $S'$ be the corresponding surface obtained from $S$. We will 
also denote the core of $\delta^{*}$ by $\delta^{*}$.
                                                                    
By Proposition 3.3, $K$ and $K'$ are $l_{\delta}$-equivalent. 
One can see that $K'$ is
$n$-hyperbolic, and it bounds an $n$-hyperbolic surface
of genus strictly less than $g$.

Obviously, there is a circle on $S'$ with $\delta^{*}$ as a sub-arc which 
bounds a disk $D$ in $S^3\setminus S'$. A surgery on $S'$ using $D$ 
changes $S'$ to $S''$ with $\partial S''=K'$, and we conclude that $S''$ 
is an $n$-hyperbolic regular Seifert surface with genus $g-1$. Thus, 
inductively, $K'$, and hence $K$, is at least $l(n)$-trivial. \qed
\medskip
\p {\bf Remark 3.3*.} 
Let us close this paragraph by remarking that with the notation 
as in Definition 3.1, if we
assume that each $[\gamma^{\epsilon}_i]$, 
$i=1,\dots,g$ and $\epsilon=+$ or $-$, lies in $\pi^{(2)}$, then
$K=\partial S$ has the trivial Alexander polynomial. Thus,
since the only
Vassiliev invariant of order two comes from the Alexander polynomial,
$K$ is {\it 2-trivial}. Moreover, examples found by 
J. Conant ([C1]) indicate that if we require
the stronger condition
that $[\gamma^{\epsilon}_i] \in \pi^{(n+1)}$ in Definition
3.1, then the lower bound in Theorem
3.2 should have linear growth on $n$.
However, Conant also has examples ([C])
demonstrating that with the weaker definition
of $n$-hyperbolicity used in this paper,
the logarithmic growth in the lower bound
of Proposition
3.3 is necessary. More precisely,
he constructs a family $\{K_n\}_{n\in \N}$
of knots such that: $K_n$ is $(2^{n+2}-4)$-hyperbolic
but,
as computer calculations indicate,
the lower bound in Proposition 3.3 should be 
$2n$. We will not pursue
an improvement of the lower bound in Theorem 3.2 here since, in the 
light of Proposition 6.2 and Conjecture 6.3, we are mainly interested
in Corollary 3.2*.
\medskip
\smallskip
\vfill
\eject
\p {\bf 3b. Nice arcs and simple commutators  }
\medskip
In this paragraph we  begin
the study of the geometric combinatorics
of arcs in good position and prove  a few auxiliary
lemmas required for the proof of Proposition 3.3.
At the same time we also describe our strategy of the proof
of 3.3.

Throughout the rest of section three, we will adapt the convention
that the endpoints of $\delta$ or of any subarc
${\tilde \delta} \subset \delta$ representing
a word in $\pi^{(m+1)}$, lie on the line $l$ associated to our fixed
projection.
Let  $W=c_1\ldots c_r$ be a word expressing $\delta^+$ as a product
of simple (quasi-)commutators of length $m+1$,
and let $p_1(S)$ be a projection of $S$, as in Lemma 2.12.
Then, each letter in $W$ is represented by a band crossing
in the projection. Now, let ${\cal C}= \{ C_1,\ldots, C_{m+1}\}$
be disjoint sets of letters obtained by applying Lemma 2.10 to the 
word $W$, so that $W$ becomes a trivial word whenever we delete letters
in a non-empty $C\in 2^{\cal C}$ from $W$ (the resulting word is denoted
by $W_C$).

Let $y$ be a free generator appearing in $W$. 
We will say that the letters $\{ y, y^{-1}\}$
constitute a {\it canceling pair}, if there is
some $C \in 2^{\cal C}$ such that the word
$W_C$ can be reduced to the identity,
in the free  group $\pi$, by
a series of deletions in which
$y$ and $y^{-1}$        
cancel with each other. 

Ideally, we would like to be able to say that
for every $C\in 2^{\cal C}$  the arc $\delta_C$
(obtained from $\delta$ by switching all crossings corresponding to
$C$) is isotopic in $S^3 \setminus S_C$ to a straight segment connecting the 
end points of $\delta$ and above the remaining diagram.
As remarked in 2.13, though, this may not always be the case.
In other words not all sets of letters $\cal C$,
that come from Lemma 2.10, will be suitable for
geometric {\it m-triviality}.
This observation leads us to the following
definition.
\medskip
\p {\bf Definition 3.4.} {\sl Let $S$, $B$ and
$\delta$ be as in the statement of 3.3
and let ${\tilde \delta}\subset {\delta}$
be a subarc that represents a word $W$ in $\pi^{(m+1)}$.
Furthermore let
$\delta^*$ be an embedded segment 
connecting the endpoints of $\tilde \delta$ and such that
$\partial {\delta^{*}} \subset l$ and
the  interior of  $\delta^*$
lies above the projection of $S$
on the projection plane.
\smallskip
\p 1) We will say that ${\tilde \delta}$
is quasi-nice if there exists a segment
$\delta^{*}$ as above and such that
either the interiors of $\delta^{*}$
and ${\tilde \delta}$ are disjoint,
or ${\tilde \delta}=\delta^{*}$ and
$\delta^{*}$ is the hook of the band $B$.
Furthermore, if
the interiors of $\delta^{*}$
and ${\tilde \delta}$ are disjoint
then $\delta^{*}$
should not separate any
set of crossings corresponding to a canceling pair
in $W$ on any of the hooks of the projection.
\smallskip
\p 2) Let $\delta$ be a
quasi-nice arc,
and let $\delta^{*}$ be as in 1).
Moreover, let $S'$ denote the surface $S\cup n(\delta^*)$,
where $n(\delta^*)$ is a flat ribbon neighborhood of
$\delta^*$. 
We will say that $\delta$ is {\it k-nice}, 
for some $k\leq m+1$, 
if there exists a collection
${\cal C}$ of $k$ disjoint sets
of band crossings on $\tilde \delta$, such that
for every non-empty $C \in  2^{\cal C}$, the loop 
$(\delta^* \cup {\tilde \delta_C})^{+}$
is homotopically trivial in
$S^3 \setminus S'_C$,
where 
$S'_C=S_C \cup n(\delta^*)$. We will say
that every $C \in  2^{\cal C}$ trivializes
$\tilde \delta$ geometrically. }
\bigskip
\centerline{\epsfxsize=3in\epsfbox{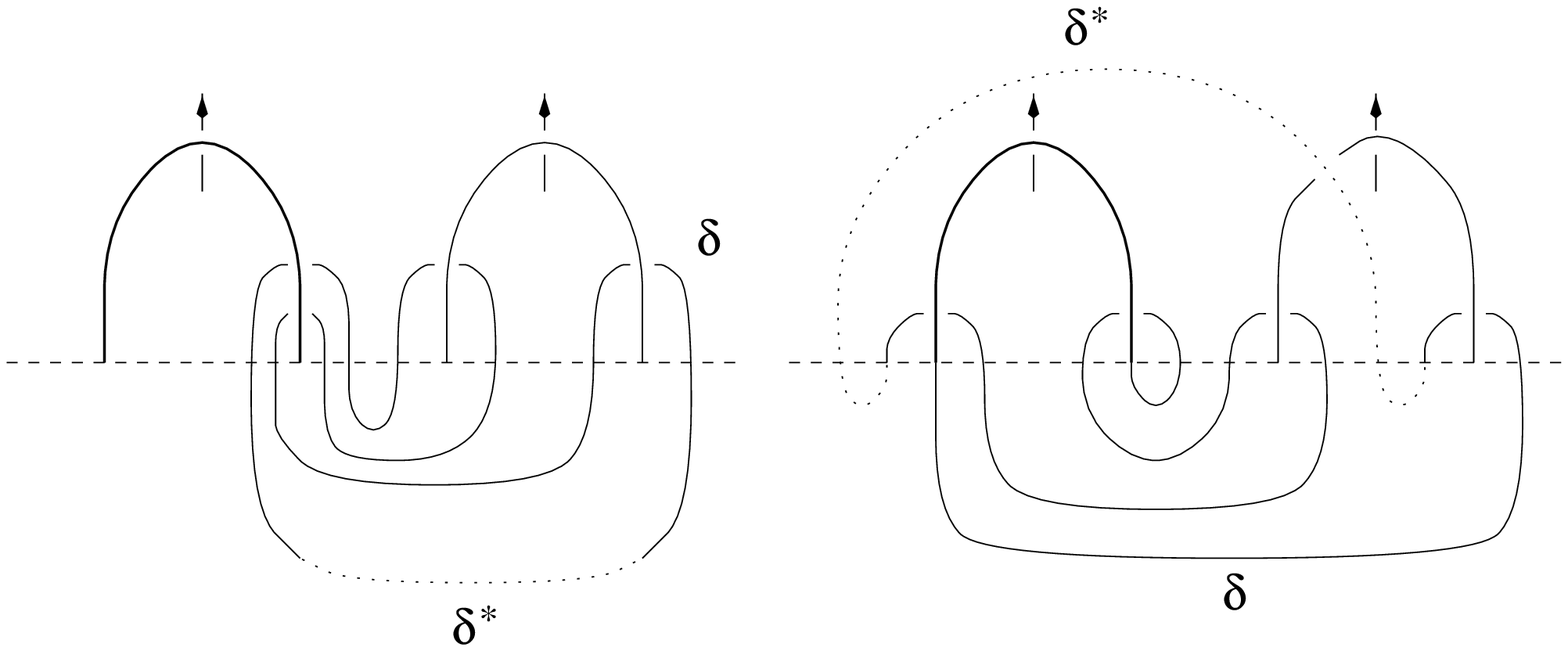} }
\medskip
\centerline {{ \bf Figure 7.} {\msmall Nice arcs representing 
simple 2-commutators}}
\medskip
Notice that the arc in the example on the left side of Figure 7
is
both {2-nice} and
quasi-nice while the one on the right side is not.
In fact, one can see that 
all embedded arcs
in good position 
representing simple 2-commutators are 2-nice.

\medskip
\p {\bf Lemma 3.5.} {\sl Let $\delta$ be a subarc of
the core of a band $B$, in a projection
of a regular surface $S$.
 Let $\delta^{*}$
be a straight segment, connecting the endpoints of $\delta$
and let $S'$ denote the surface $S\cup n(\delta^{*})$.
Suppose that the loop $(\delta^{*} \cup \delta)^+$
is homotopically trivial in
$S^3 \setminus S'$. Then $n(\delta)$
can be isotoped onto $n(\delta^{*})$ in
$S^3 \setminus S$ relative to the endpoints.}
\medskip

\p {\it Proof.}  Since $\delta^{*} \cup \delta \subset S'$
is an embedded loop, by Dehn's Lemma (see for example [He] or [Ro])
 we conclude
that it bounds an embedded disc
in $S^3 \setminus S'$.
Then the claim follows easily. \qed

\medskip
\p {\bf Corollary 3.6.} {\sl Let $\gamma$ and $\delta$
as in the statement of Proposition 3.3.
Assume that
$\delta$ is an $q_{\delta}$-nice arc.
Then the conclusion of the Proposition is true for $\delta$.}
\medskip

\p {\it Proof.} It follows immediately from
Definitions 1.2 and 3.4 and Lemma 3.5. \qed
\medskip
For the rest of this paragraph we will
focus on projections  of arcs, in {\it good position},
that represent simple quasi-commutators.
We will analyze the geometric combinatorics of such projections.
This analysis will be crucial, in the next paragraphs,
in showing
that an arc $\delta$ as in Proposition 3.3 
is $q_{\delta}$-{nice}.

\smallskip
Let $\delta_1$ be a subarc
of $\delta$ presenting a simple quasi-commutator
of length $m$, say $c$. Moreover, let 
$\delta_2$ be another subarc of $\delta$ presenting a simple
quasi-commutator equivalent to $c$ or $c^{-1}$.
We may change the orientation of $\delta_2$ if necessary
so that it presents a simple quasi-commutator equivalent to $c$. Then we may
speak of the initial (resp. terminal) point $p_{1,2}$
(resp. $q_{1,2}$) of $\delta_{1,2}$; recall these points all lie on the 
line $l$. 
\medskip
\p {\bf Definition 3.7.} {\sl Let ${\hat \delta_1}$
(resp. ${\hat \delta_2}$) be the segment on $l$ going from $p_1$ to $p_2$
(resp. $q_1$ to $q_2$). We say that $\delta_1$ and $\delta_2$
are {\it parallel} if the following are true:
i) At most one hook has its end points on $\hat\delta_1$ or $\hat\delta_2$
and both of its end points can be on only one of $\hat\delta_{1,2}$;
ii) If a hook has exactly one point on some $\hat\delta_j$,
say on $\hat\delta_1$, then $\hat\delta_1$
dosen't intersect the interior of
$\delta_{1,2}$.
iii) We have either $\hat\delta_1\cap\hat\delta_2=\emptyset$
or $\hat\delta_1\subset \hat\delta_2$;
iv) If $\hat\delta_{1,2}$ are drawn disjoint and above the surface $S$,
the diagram $\delta_1\cup\hat\delta_1\cup\delta_2\cup \hat\delta_2$
can be changed to an embedding by type II Reidermeister moves.
}
\medskip 
The reader may use Figure 8 to understand Definition 3.7. It should not be 
hard to locate the arcs $\hat\delta_{1,2}$ in each case in Figure
8. In 
the first two pictures, 
the straight arcs $\hat\delta_{1,2}$ have no crossings with
$\delta_{1,2}$. Crossings between $\hat\delta_{1,2}$ and $\delta_{1,2}$
removable by type II Reidermeister moves are allowed to accommodate the 
modification of $\delta_{1,2}$ in Lemma 2.12.
For example, in the last two pictures of Figure
8 one of $\hat\delta_{1,2}\subset l$
intersects both of $\delta_{1,2}$.

\bigskip
\centerline{\epsfxsize=3in\epsfbox{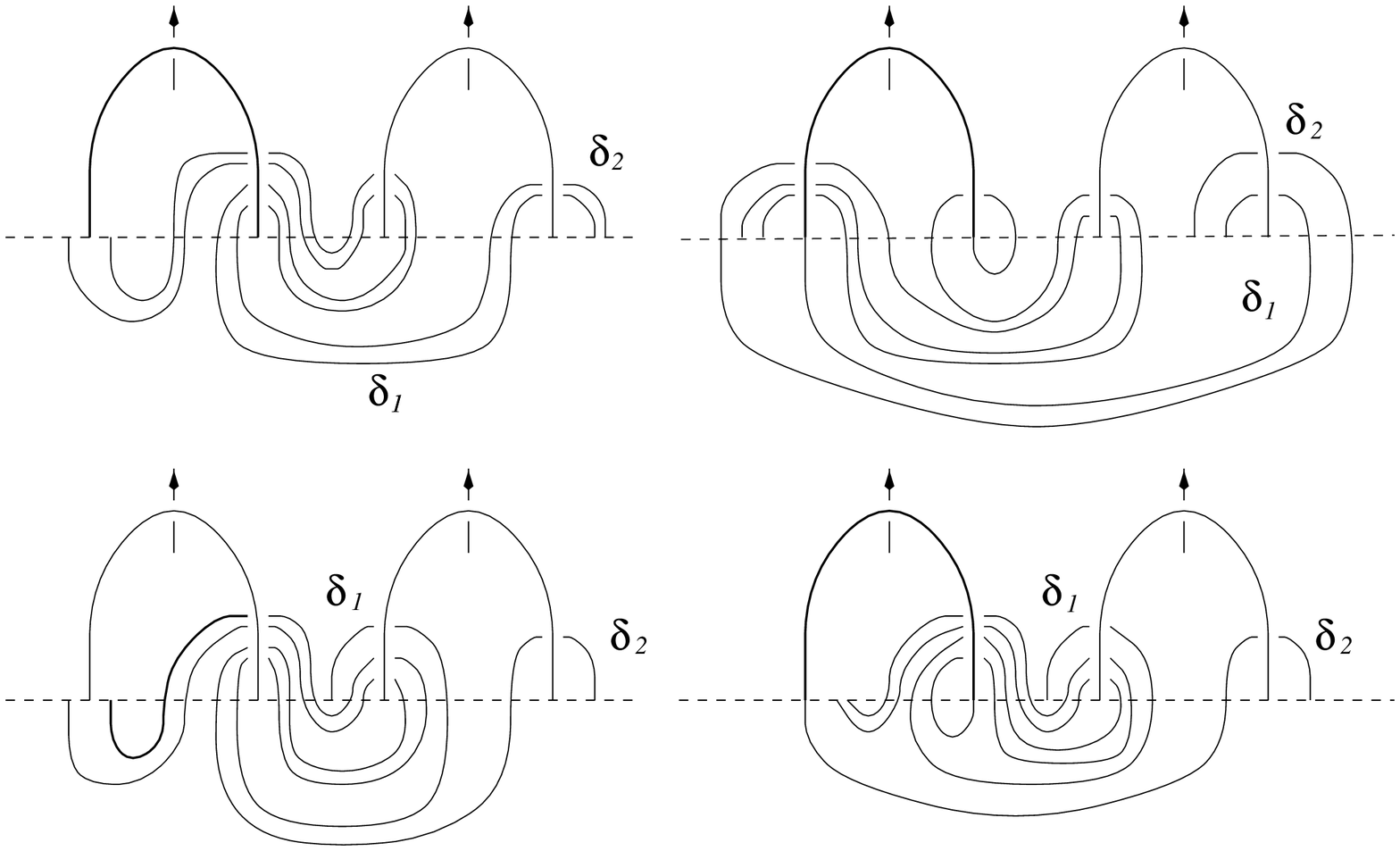} }
\medskip
\centerline {{ \bf Figure 8.} {\msmall Various kinds of parallel arcs}}
\medskip

\medskip
\p {\bf Lemma 3.8.} {\sl Assume that the setting is as in the statement of 
Proposition 3.3. Let $c_1$ and $c_2$ be equivalent simple quasi-commutators 
presented by sub-arcs $\delta_{1,2}$ of $\delta$ respectively, and let $y$ be 
one of the free generators associated to the hooks of our fixed projection.
Moreover, assume that $\delta_{1,2}$ are parts of a subarc $\zeta$ of 
$\delta$ presenting a simple quasi-commutator 
$W=c_1yc_2^{-1}y^{-1}$. Then $\delta_1$ and $\delta_2$ are parallel.}
\medskip

\p {\it Proof.} By abusing  the notation, we denote 
$\delta_1=\tau_1 x \mu_1 x^{-1}$ and $\delta_2=\tau_2 x \mu_2x^{-1}$ where 
$\tau_1,\tau_2,\mu_1^{-1},\mu_2^{-1}$ all present equivalent 
simple quasi-commutators. 
Furthermore,
$\zeta=\delta_1y\delta_2^{-1}y^{-1}$. 
For a subarc $\nu$, up to symmetries, there are four possible ways 
for both of its endpoints 
to reach a certain $z$-hook so that $z\nu z^{-1}$ is presented by an arc 
in good position. 
See Figure 9, where the arc $\nu$ may run through the 
$z$-hook. We will call the pair of undercrossings $\{z,z^{-1}\}$ a {\it
canceling pair}.
Now let us consider the relative positions of $\tau_1$, $\tau_2$, $\mu_1$ and
$\mu_2$. Inductively, $\tau_i$ and $\mu_i^{-1}$ are parallel, for $i=1,2$.
If $x\mu_1x^{-1}$ is of type (I) in Figure 9, since $\tau_1$ and
$\mu_1^{-1}$ 
are parallel, $\tau_1$ has to go the way indicated in Figure 10 (a).
\bigskip
\centerline{\epsfxsize=2in\epsfbox{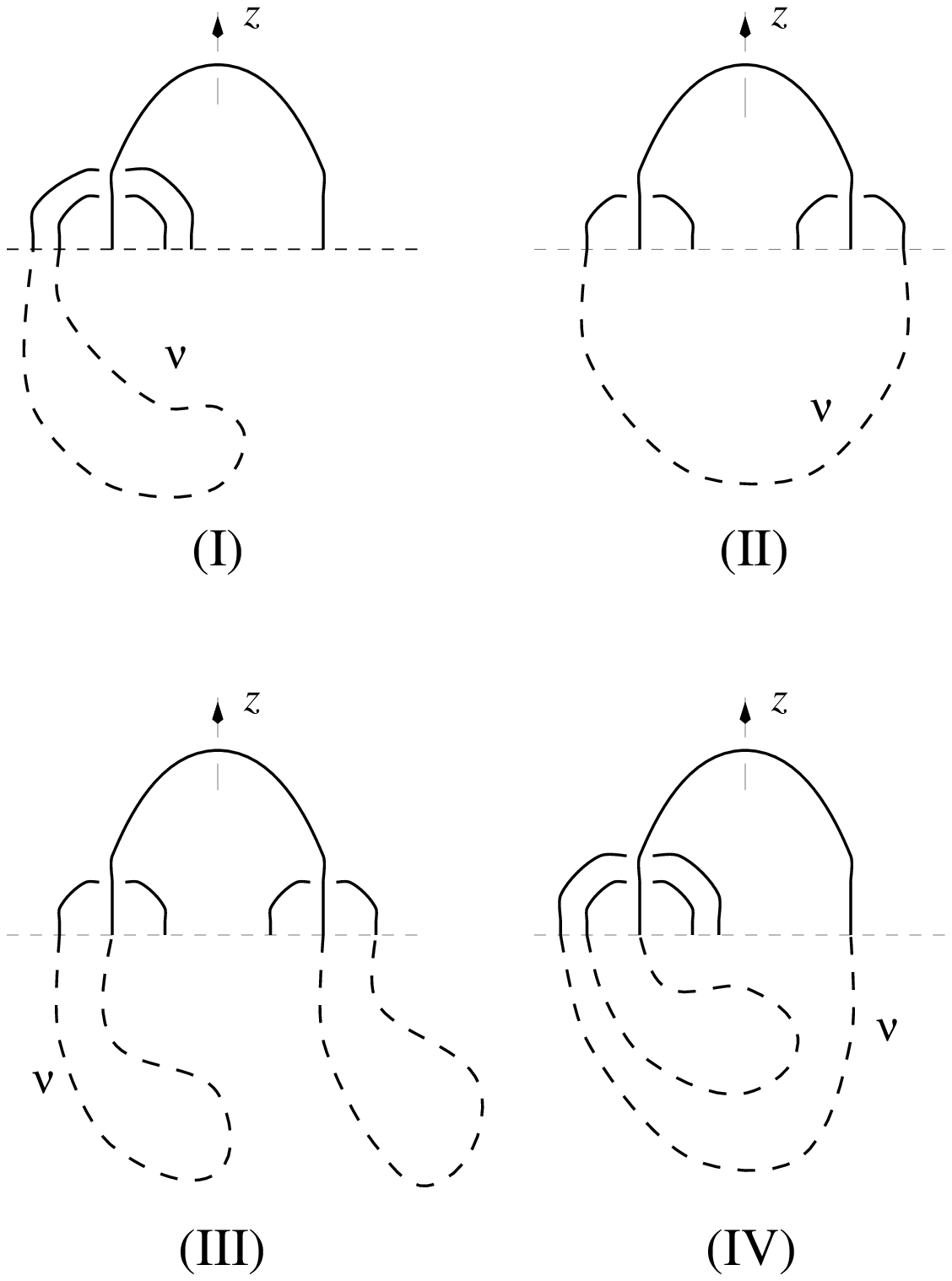} }
\medskip
\centerline {{ \bf Figure 9.} {\msmall Types of arcs presenting $z\nu z^{-1}$}}
\medskip
If 
$x\mu_2x^{-1}$ is also of type (I), there are two cases to consider. One case 
is to have the canceling pairs $\{x,x^{-1}\}$ in $x\mu_1x^{-1}$ and 
$x\mu_2x^{-1}$ both going underneath the $x$-hook at the left side, and the 
other case is to have them going underneath the $x$-hook at different sides.
In the first case, in order to read the same word from $\tau_1$ and $\tau_2$ 
as well as from $\mu_1$ and $\mu_2$, $\tau_1x\mu_1x^{-1}$ and 
$\tau_2x\mu_2x^{-1}$ has to fit like in Figure 10 (b).
This implies that 
$\delta_1$ and $\delta_2$ are parallel. 
In the second case 
(see Figure 10 (c)), in order that $\delta_i$ be parts of the arc 
$\zeta=\delta_1y\delta_2^{-1}y^{-1}$, they have to go to reach the same 
$y$-hook.

\medskip
\centerline{\epsfxsize=2.7in\epsfbox{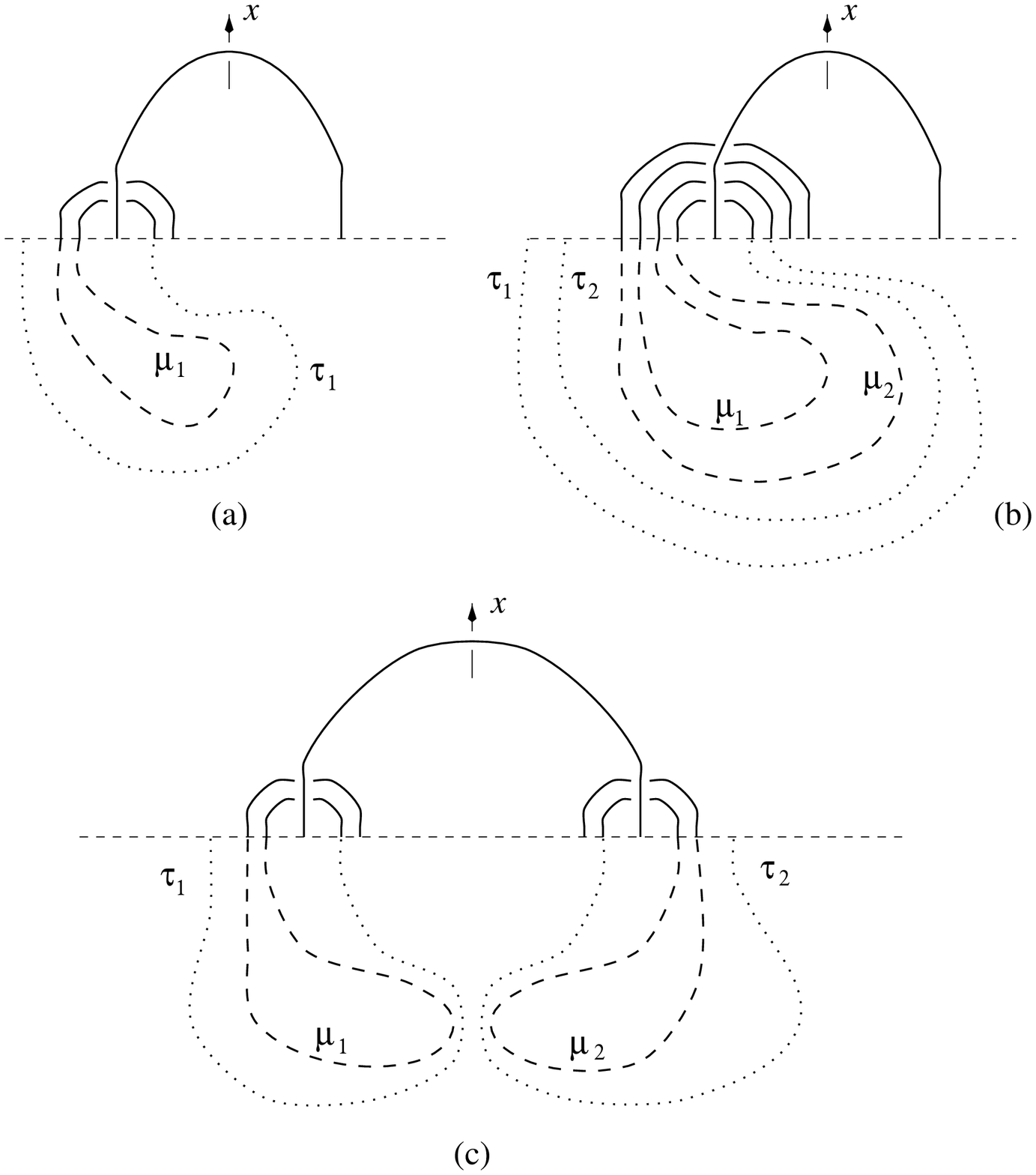} }
\medskip
\centerline {{ \bf Figure 10.} {\msmall The case when both $x\mu_1x^{-1}$ and 
$x\mu_2x^{-1}$ are of type (I)}}
\medskip
But then we will not
be able to
read the same word through $\tau_1$ and $\tau_2$. This shows that if 
$x\mu_1x^{-1}$ and $x\mu_2x^{-1}$ are both of type (I), $\delta_1$ and
$\delta_2$ are parallel.
There are many other cases which can be checked one by one in the same way as
in Figure 10.
The details are left to  the  patient reader. \qed
\medskip
So far we have been considering the
projection of our surface on a plane  $P$
inside ${\R}^3= P \times {\R}$.
To continue, let us  pass to the compactifications of ${\R}^3$
and $P$. We obtain a 2-sphere $S^2_{P}$ inside $S^3$,
and assume that our projection in Proposition
3.3 lies on $S^2_{P}$. We may identify the image of $l$
with the equator of $S^2_{P}$,
and the images of $H_+$ and $H_-$
with the upper and lower hemisphere. 
We will interchange between $P$ and $S_P^2$
whenever convenient.

\medskip
\p {\bf Digression 3.9.}\ Let $\delta$, $B$ be as in the statement 
of Proposition 3.3
and let $x_0$ denote the free generator of $\pi_1(S^3\setminus S)$ 
corresponding to $B$.
Suppose
$\delta_1$ and $\delta_2$ are parallel subarcs of
$\delta$. and let $\hat\delta_{1,2}$
be as in Definition 3.7.
We further assume that the 
crossings between
$\hat\delta_{1,2}$ and $\delta_{1,2}$ have been removed by isotopy. 
Let $y$ be a free generator 
of $\pi=\pi_1(S^3 \setminus S)$. We assume that both $\delta_1$ and $\delta_2$
are sub-arcs of an arc $\zeta \subset \delta$ presenting $[c^{\pm 1},\ 
y^{\pm 1}]$ (recall that $\delta_{1,2}$ present $c^{\pm1}$). 
Then $\zeta$ is 
a union $\delta_1\cup\tau_1\cup \delta_2\cup\tau_2$, where $\tau_{1,2}$ are
segments each going once underneath the $y$-hook. 
One point of $\zeta$ is the same as one endpoint of one
of $\delta_{1,2}$, say $\delta_2$.
Let ${\bar \delta}:=\delta \setminus (\delta_1 \cup \delta_2)$.
By the properties of good position we see
that in order for one of $\delta_{1,2}$, say $\delta_1$,
not to be embedded 
on the projection plane it must run through the hook part of $B$,
and the word representing $\delta_1$ must involve $x_0$.
Moreover, good position imposes a set
of restrictions on the relative positions of
of $\delta_{1,2}$
and the various subarcs of ${\bar \delta}$.
Below we summarize the main 
features of
the relative positions of $\delta_{1,2}$
and the various subarcs of ${\bar \delta}$; these
features will be useful
to us in the rest of the paper.
We will mainly focus on the case
that $\delta_{1,2}$ are embedded; the case
of non-embedded arcs is briefly discussed in part 
b)
of this Digession.
\smallskip
\p a) Suppose that $\delta_{1,2}$ are embedded
on the projection plane $P$.
Then the loop  
$\delta_1 \cup {\hat \delta_2} \cup \delta_2 \cup {\hat \delta_1}$ 
separates $S_P^2$ into two discs,
$D_1$ and $D_2$.
 The intersections $D_{1,2}\cap{\bar \delta}$ consist of finitely
many arcs. With the exception
of at most one 
these arcs are embedded.
One can see (see the two pictures on the left side
of Figure 8) that the interiors of
$\tau_{1,2}$
are 
disjoint from that of exactly one of $D_{1,2}$, say 
$D_1$, and they lie in the interior of the other.
We will call $D_1$ (resp. $D_2$ ) the {\it finite}
(resp. {\it infinite}) disc corresponding to the pair
$\delta_{1,2}$.
Using the properties of 
good position one can see that for each component
$\theta$ of $D_1 \cap {\bar \delta}$ one of the  following is true:

\p $(a_1)$  Both the endpoints of $\theta $ lie on
${\hat \delta_2}$ and $\theta $ 
can be pushed in the infinite disc $D_2$ after isotopy,
or it represents a word $w$,
such that the following is true:
None of the letters appearing
in the reduced form of $w$ appears
in the underlying commutator of
$c$. Moreover for each
free generator $x$ appearing
in $w$, 
$C^{\pm 1}$ contains  inserted pairs $x^{\pm 1} x^{\mp 1}$.
To see these claims, first notice that
if one of the generators, say $z$,
appears in the underlying commutator
of $c$ then the intersection
of $D_1$ and the $z$-hook
consists of two (not necessarily disjoint)
arcs,
say $\theta_{1,2}$,
such that on point of
$\partial (\theta_{1,2})$
is on $\delta_1$ and the other on $\delta_2$.
Moreover, both the endpoints of the
$z$-hook lie outside $D_1$ in the infinite disc.
Now a subarc of ${\bar \delta}$
in $D_1$ has the choice of
either hooking with
$z$ in exactly the same fashion as
$\delta_{1,2}$, or
``push"  $\theta_{1,2}$ by
a finger move as indicated in Figure 4,
and hook with some $x \neq z^{\pm 1}$.
In order for the second possibility to occur,
at least one of the endpoints of the $x$-hook
must lie inside $D$; by our discussion
above this will not happen if
$x$ has already appeared in the underlying
commutator of $c$. The rest
of the claim follows from the fact that
the ``top" of the $x$-hook must lie
outside $D_1$
\medskip
\p $(a_2)$ One endpoint of $\theta$ is
on
${\hat \delta_1}$ and the other on ${\hat \delta_2}$.
Moreover,
$\theta$ either
represents $c^{\pm 1}$ or the trivial word
or it represents a $w$, as in case
$(a_1)$ above.

\medskip
\p $(a_3)$ $\theta$ is a subarc of the hook
corresponding to the band $B$,  and it has one endpoint on
$\delta_1$ and the other on $\delta_2$
or one point on
$\delta_i$ and the other
on ${\hat \delta_j}$ ($i,j=1, 2$). Furthermore,
if $\theta$ has one endpoint on $\delta_i$ and the other
on ${\hat \delta_j}$ then i) the underlying commutator of $c$
does not involve $x^{\pm 1}_0$;
and ii) both endpoints
of the $x_0$-hook lie inside $D_1$.
\medskip
\p $(a_4)$ One of $\delta_{1, 2}$, say $\delta_1$, runs through
the hook part of $B$ and $\theta$
has one point on
$\delta_1$ and the other
on ${\hat \delta_j}$ ($i,j=1, 2$). Moreover,
we have:  
i)The arc
$\theta$ represents $w x_0^{\pm 1}$
where either $w=c^{\pm 1}$ or
none of the letters appearing
in there reduced form of $w$ appears
in the underlying commutator of
$c$; ii) the underlying commutator of
$c$ does not involve $x_0^{\pm 1}$; iii) if
$e$ is a
simple quasi-commutator represented
by a subarc ${\tilde \delta}$ such that
${\delta_{1,2}} \subset {\tilde \delta}$,
then the underlying commutator
of $e$ does not involve $x_0^{\pm 1}$ (see
also Lemma 3.12 a)).

\medskip
\p b) Recall that $x_0$ is the free generator of $\pi$
corresponding to $B$.
Suppose that $\delta_1$ is non-embedded.
Then $\delta_1$ to run through the hook of $B$
and the word representing $\delta_1$
must involve $x_0$.
\smallskip
\p $(b_1)$ It follows from the properties
of good position that
any subarc of $\theta \subset {\delta}$ that has its
endpoints on different ${\hat {\delta_i}}$ has 
to represent $c^{\pm 1}$.
\medskip
\p {\bf Lemma 3.10.} {\sl Let
the setting be again as in the statement of 
Proposition 3.3, and let $\delta_1$ be a subarc
of $\delta$ representing a simple (quasi-)commutator.
Moreover let $h_0$ denote the hook
part of the band $B$.
We can connect
the endpoints of $\delta_1$
by an arc $\delta_1^*$,
which is embedded on the projection plane
and such that: i) $\delta_1^*$ lies on the top
of the projection $p_1(S)$; 
ii) the boundary $\partial (\delta_1^*)$ lies on the line $l$;
and
iii) either $\delta_1^*=h_0$ 
or  the interiors
of $\delta_1^*$ and
$\delta_1$ are disjoint
and
$\delta_1^*$ goes over at most one hook at most once.}
\medskip

\p {\it Proof.} Suppose that $\delta_1$ represents
$W=[c, y^{\pm 1}]$ and let $\delta_1^{1,2}$
be the subarcs of $\delta_1$
representing $c^{\pm 1}$ in $W$. By Lemma 3.8
$\delta_1^{1,2}$ are parallel; let ${\hat \delta_1}^{1,2}$
be the arcs of Definition 3.7 connecting the endpoints
of $\delta_1^{1,2}$. Recall that there is at most one hook,
say corresponding to a generator $z$, that
can have its endpoints on $\delta_1^{1,2}$. If
$z=y$ then , using good position, we see that
there is an arc $\delta_1^*$ as claimed above
such that either $\delta_1^*=h_0$
or it 
intersects at most the $y$-hook in at most 
one point. If $z\neq y$
then we can find an arc $\alpha$
satisfying i) and ii) above and such that
either $\alpha=h_0$ or $\alpha$ intersects
the $y$-hook in at most 
one point and  the intersections of $\alpha$
with the other hooks can be removed
by isotopying $\alpha$, relatively 
its endpoints. Thus the existence of
$\delta_1^*$ follows again. \qed

\medskip
\p {\bf Definition 3.11.} 
{\sl An arc $\delta_1$ representing a simple (quasi-)commutator
will be called {\it good} if the arc $\delta^*$
of Lemma 3.10, connecting the endpoints
of the $\delta_1$, doesn't separate
any canceling pair of crossings in $\delta_1$.}
\medskip
The reader can see that the arc in
the picture of the left side
of Figure 7 is {good} while the one on the right 
is not {good}.
\medskip
\p {\bf The outline of the proof of 3.3.}
In 3.4 we defined the notions of quasi-niceness and
$k$-niceness. By definition, a $k$-nice arc is quasi-nice.
We will, in fact, show that the two notions are equivalent.
More precisely, we show in Lemma 3.21 that a
quasi-nice arc $\delta$ is $q_{\delta}$-nice.
This, in turn, implies Proposition 3.3.
To see this last claim,
 notice that the arc $\delta$ in the statement
of 3.3 is quasi-nice. Indeed, 
since the arc $\eta$
in the statement of 3.3 represents the identity
in $\pi$, good
position and the convention about the 
endpoints of $\delta$ made
in the beginning of 3b assure the following:
Either the interior of $\eta$
lies below $l$ (and above
the projection of $S\setminus n(\eta)$)
and it is disjoint from that of
$\delta$ or $\eta$ is the hook
part of $B$. In both cases 
we choose $\delta^{*}=\eta$.

To continue, notice that a
good arc is by definition quasi-nice.
The notion of a {\it good } arc
is useful
in organizing and studying the 
various simple quasi-commutator pieces
of the arc $\delta$ in 3.3. 
In Lemma 3.15 we show that if an arc $\tilde \delta$ 
is good then it is $q_{\tilde \delta}$-nice
and in Lemma 3.19 we show
that if $\tilde \delta$ is a product of good arcs
then it is $q_{\tilde \delta}$-nice. In both cases 
we exploit good position to estimate the number 
of ``bad" crossings along ${\tilde \delta}$,
that are suitable for algebraic triviality but may 
fail for geometric triviality. 
All these are done in 3c and 3d.

In 3e we begin with the observation
that if $\tilde \delta$ is a
product of arcs $\theta_1,...\theta_s$
such that $\theta_i$ is $q_{\theta_i}$-nice
then $\tilde \delta$ is $q_{\tilde \delta}$-nice (see Lemma 3.20).
Finally, Lemma 3.21 is proven by induction on the number of
``bad" subarcs that $\delta$ contains.
\vfill
\eject
\medskip
\p {\bf 3c. Bad sets
and good arcs } 
\medskip
In this paragraph we continue our
study of arcs in good position that represent
simple quasi-commutators.
Our goal, is to show that
a good arc $\delta$ representing a
simple quasi-commutator is $q_{\delta}$-nice (see
Lemma 3.15).

Let $W=[\ldots [[y_{1},\  y_{2}],\  y_{3}], \ldots,  y_{m+1}]$ 
be a simple (quasi-)
commutator represented by an
arc
$\delta$ of the band $B$ in good position.
Suppose that  the subarc of $\delta$ representing
$W_1=[\ldots [[y_{1},\  y_{2}],\  y_{3}], \ldots,  y_i]$,
for some $i=1, \ldots, m+1$,
runs through the hook part of $B$, at the stage that
it realizes the crossings
corresponding to $y_i$. The canceling
pair corresponding to $\{y_i,\  y_i^{-1}\}$
will be called {\it the special canceling pair}
of $W$. 
\smallskip
\p {\it Note:} Let $x_0$ be the free generator
of $\pi$ corresponding to the band $B$.
Notice that if the special
canceling pair $\{y_i,\  y_i^{-1}\}$
is of type (III) or (IV)
then we must have $y_i=x_0^{\pm1}$;
otherwise we must have $y_i\neq x_0^{\pm1}$.
Moreover, it follows by good position and Lemma
3.8 that if $\{y_i,\  y_i^{-1}\}$
is of type (I) or (II), then the arc
representing $W_1=[\ldots [[y_{1},\  y_{2}],\  y_{3}], \ldots,  y_{i-1}]$
is embedded; that is $y_j \neq x_0^{\pm 1}$ for $j=1, \ldots, i-1$.
\medskip
\p {\bf Lemma 3.12.} {\sl Let 
$W=[\ldots [[y_{1},\  y_{2}],\  y_{3}], \ldots,  y_{m+1}]$ 
be a simple (quasi-)
commutator represented by an arc
$\delta$ of the band $B$ in good position,
and let
$x_0$ be the free
generator of $\pi$ corresponding to the hook of $B$.

\p a) Suppose that, for some $i=1, \ldots, m+1$,
one of the canceling pairs
$\{y_i,\  y_i^{-1}\}$
is the special canceling pair.
Then, we have $y_j \neq x_0^{\pm 1}$
for all $i<j\leq m+1$.

\p b) Let $z$ be any free generator 
of $\pi$
corresponding to one of the hooks of our projection.
Then, at most two successive 
$y_{i}$'s can be equal to $z^{\pm 1}$. }
\medskip
\p {\it Proof.} a) For $j>i$
let $c= [\ldots [[y_{1},\  y_{2}],\  y_{3}], \ldots,  y_{j-1}]$ and let
$\delta_{1,2}$
be the arcs representing  $c^{\pm 1}$
in $[c, \ y_j]$. By 3.8, $\delta_{1,2}$
are parallel. Let
$\hat \delta_{1,2}$
be arcs satisfying Definition 3.7.
Notice that the $x_0$-hook can not have just
one of its endpoints on $\hat \delta_{1,2}$.
For, if the $x_0$-hook had one endpoint on, say,
$\hat \delta_{1}$, then
$\hat \delta_{1}$ would intersect
the interior of $\delta_{1,2}$.
Now easy drawings will convince us that
we can not form $W_1=[c, \ x_0^{\pm 1}]$
without allowing the arc representing
it to have self intersections below the line $l$
associated to our projection.
But this would
violate the requirements
of good position.
\medskip
\p b)
By symmetry we may assume  that
$W=[\ldots, [c^{\pm 1},\ z^{\pm 1}], y_i,  \ldots, y_{m+1}]$,
where $c$ is a simple (quasi-)commutator
of length $<m$. Let $\tilde \delta$
be the subarc of $\delta$ representing $[c^{\pm 1},\ z^{\pm 1}]$,
and
let $\delta_{1,2}$ 
be the subarcs of $\tilde \delta$ representing $c^{\pm 1}$.
By Lemma 3.8, $\delta_{1,2}$ are parallel.
Let ${\hat \delta}_{1,2}$ be as in Definition 3.7.

Without loss of generality we may assume that
$z^{\pm 1}$ has already appeared in $c$. Thus the intersection
$\delta_{1,2} \cap H_{+}$ is a
collection of disjoint arcs $\{A_i \}$, each passing once under the
$z$-hook, and with their endpoints on the line $l$.
Let $A_{1,2}$ denote the innermost of the $A_i$'s
corresponding to the two endpoints of the $z$-hook.
Let  $\alpha_{1,2}$ denote the
segments of $l$ connecting the endpoints
of $A_{1,2}$, respectively.
Observe
that 
both of the endpoints of at least one 
of ${\hat \delta}_{1,2}$ must lie
on $\alpha_{1}$ or $\alpha_{2}$. There are
three possibilities:
\smallskip
(i) The endpoints of both ${\hat \delta}_{1,2}$
lie on the same $\alpha_{1,2}$, say on $\alpha_{1}$; 

(ii) The endpoints of  ${\hat \delta}_{1,2}$
lie on different $\alpha_{1,2}$;

(iii) The endpoints of one of ${\hat \delta}_{1,2}$
lie outside the endpoints of
$\alpha_{1,2}$.

Suppose we are in (i) above: Notice that both
the endpoints of the arc $\tilde \delta$
must also lie on $\alpha_{1}$. By Definition 3.7 we see
that both of the endpoints of any arc parallel to
$\tilde \delta$ must also lie on $\alpha_{1}$.
There are two possibilities for the relative 
positions of ${\hat \delta}_{1,2}$;
namely $\hat\delta_1\cap\hat\delta_2=\emptyset$
or $\hat\delta_1\subset \hat\delta_2$.
Using Digression 3.9 we can see that in both cases
we must have $y_i \neq z^{\pm 1} $. 

If the last
letter in $c^{\pm 1}$ was not $z^{\pm 1}$ then we are 
done. Suppose, now, that $c=[d, \ z^{\pm 1 }]$.
Then, $A_1$ must be part of one of $\delta_{1,2}$
and one endpoint
of one of ${\hat \delta}_{1,2}$
must be an endpoint of $\alpha_1$. 
Moreover, observe that the last canceling pair in $c$
can not be of type (I).
By part a)
this last canceling pair in $c$ must be of type (II), 
and the only way
that the last letter in $d$ can be $z^{\pm 1}$
is to have $d=[e, \ z^{\pm 1 }]$ where
the last
canceling pair is of type (I) or (II). Now easy drawings, 
using Digression 3.9, will convince us
that this is not possible.

We now proceed with case (ii) above: Observe that
the last letter in $c$ can
not be $z^{\pm 1}$; otherwise 
the endpoints of each of $\delta_{1,2}$ would be
on the same $\alpha_{1,2}$. We first form $[c^{\pm 1},\ z^{\pm 1}]$.
The endpoints of the arc $\tilde \delta$
representing $[c^{\pm 1},\ z^{\pm 1}]$ are now on the same 
$\alpha_{1,2}$, say on $\alpha_{1}$. Thus both of the endpoints of 
any arc parallel to
$\tilde \delta$ must also lie on $\alpha_{1}$. Now we are in
the situation of (i) above and the conclusion follows.

Finally, assume we are in case (iii) above: Notice that 
both endpoints
of the arc $\tilde \delta$
representing $[c^{\pm 1},\ z^{\pm 1}]$ lie
outside $\alpha_{1,2}$. Thus $y_i \neq z^{\pm 1}$
or we are in the situation of (i) above.
Using an argument similar to these of (i), (ii)
we can see that we can't have the
two last letters of $c^{\pm 1}$ being equal
to $z^{\pm 1}$. \qed
\medskip
In order to continue
we need some notation and terminology.
We will write $W= [y_{1},\  y_{2},\  y_{3}, \ldots,  
y_{m+1}]$ to denote the 
simple (quasi-)commutator 

$$W=[\ldots [[y_{1},\  y_{2}],\  y_{3}], \ldots,  y_{m+1}].$$
Let $C_1,\ldots, C_{m+1}$ be  the sets of letters 
of Lemma 2.10, for $W$. 
Recall that for every
$i=1, \ldots, m+1$ the only letter appearing in $C_i$
is $y_{i}^{\pm 1}$. 
\medskip
\p {\bf Definition 3.13.} {\sl
We will say that the set $C_i$
is {\sl bad}
if there is some $j\neq i$ such that
i) we have
$y_j=y_i=y$, for some free generator $y$;
and ii)
the crossings on the
$y$-hook,  corresponding 
to a canceling
pair $\{ y_j, \ y_j^{-1} \}$ in $C_j$,
are separated by crossings in $C_i$.}
\smallskip
The problem with a bad $C_i$
is that changing the crossings in $C_i$
may not trivialize the arc
$\delta $ geometrically.
\medskip
For $i=2, \ldots, m+1$, 
let $c_i= [y_{1}, \ \ldots,  
y_{i-1}]$ and let $\delta_{1,2}$
be the parallel arcs representing
$c_i^{\pm 1}$  in $[c_i, y_i]$.
Let $\bar \delta=
\delta\setminus (\delta_1 \cup \delta_2)$.
We will say that
the  canceling pair $\{ y_i, y_i^{-1}\}$
is {\sl admissible} if
it is of type (I) or (II).
\medskip
\p {\it Note:} If $W$ doesn't involve the generator
of $\pi$ corresponding to the hook of $B$
then every canceling pair
in $W$ is {\sl admissible}.
\medskip
\p {\bf Lemma 3.14.} {\sl 
\p a) Let $W=[y_{1},\  y_{2},\  y_{3}, \ldots,  
y_{m+1}]$ be a
simple quasi-commutator represented
by an arc $\delta$ in good position and
let $z$
be a free generator.
Also, let $C_1,\ldots, C_{m+1}$
be sets of letters as above.
Suppose that $C_i$
is bad and let $\{ y_j, \ y_j^{-1} \}$
be a canceling 
pair in $C_j$, whose
crossings on the
$z$-hook
are separated by crossings in $C_i$.
Suppose, moreover, that
the pair $\{ y_i, \ y_i^{-1} \}$
is admissible.
Then, with at most one exception,
we have $j=i-1$ or $j=i+1$.
\medskip
\p b) Let $w(z)$
be the number of the $y_i$'s in $W$ that are equal to
$z^{\pm 1}$. There can be
at most ${\displaystyle {\left[{{w(z)}\over 2}\right]
+1}}$
 bad sets
involving $z^{\pm 1}$. 
\medskip
\p c)For every $j=2, \ldots, m+1$,
at least one of $[y_1,\ \ldots,  
y_{j-1}]$ and $[y_1,\ \ldots,  
y_j ]$ is represented by a good arc.}
\medskip
\p {\it Proof.} a) Let $c= [y_{1}, \ \ldots,  
y_{i-1}]$, let $\delta_{1,2}$
be the parallel arcs representing
$c^{\pm 1}$  in $[c, y_i]$ and let
${\hat \delta_{1,2}}$ be the arcs of Definition 3.7.
Let $C_z$ denote the canceling pair
corresponding to $y_i^{\pm 1}$
in $[c,\ y_i]$.
Moreover, let $\bar \delta=
\delta\setminus (\delta_1 \cup \delta_2)$
and let 
${\bar \delta}_c$ denote the union of arcs
in 
${\bar \delta}$ such
that i) each has one endpoint on ${\hat \delta}_1$
and one on ${\hat \delta_2}$ and 
ii) they do not represent copies of $c^{\pm 1}$
in $W$. By Lemma 3.8, and Digression 3.9,
it follows that ${\bar \delta}_c = \emptyset$.

Without loss of generality we may assume that
$j>i$. Also, we may, and will,
assume that $y_{i+1} \neq z^{\pm 1}$.
\smallskip
First suppose that $C_z$ is of type (II):
By Lemma 3.12(a), it follows that if
one of $\delta_{1,2}$ has runned through
the hook part of $B$
then ${\bar \delta} \cap \delta_{1,2}= \emptyset$..
Thus the possibility discussed in $(a_4)$
of Digression 3.9 doesn't occur.
Now, by Lemma 3.8, it follows
that in order for the crossings corresponding to 
$C_z$ to separate
crossings corresponding to a
later appearance of $z^{\pm 1}$, we must have 
i) $y_j= z^{\pm 1}$ realized by a
canceling pair of type (II) and ii)
the crossings on the $z$-hook
corresponding to
$y_j$ and $y_j^{-1}$ lie below (closer to endpoints
of the hook) these representing $y_i$ and $y_i^{-1}$.
But a moment's thought,
using 3.9 and
the assumptions
made above, will convince us that in order for 
this to happen we must have ${\bar \delta}_c \neq \emptyset$; 
which is impossible. 
\smallskip
Suppose $C_z$ is of type (I):
Up to symmetries, the configuration for 
the arc $\tilde \delta$, representing $[c,\ y_i]$,
is indicated in  Figure 10(b). The
details in this case are similar to the previous case
except that
now the two crossings 
corresponding to $\{ y_i, \ y_i^{-1}\}$
occur on the same side of the $z$-hook
and we have the following
possibility: Suppose
that $c$ does not contain any type (II)
canceling pairs in $z$.
Then, we may have a type (II)
canceling pair $\{ y_j, \ y_j^{\pm 1} \}$,
for some $j>i$, such that
the crossings in
$C_z$ separate
crossings corresponding 
$y_j$ and
${\bar \delta_c}= \emptyset$.
This corresponds to the exceptional case mentioned in the statement
of the lemma. In this case, the arc
representing 
$[y_1,  \ldots,  y_i,\ldots, y_j]$
can be seen to be a {\sl good} arc.
\medskip
\centerline{\epsfysize=1.6in\epsfbox{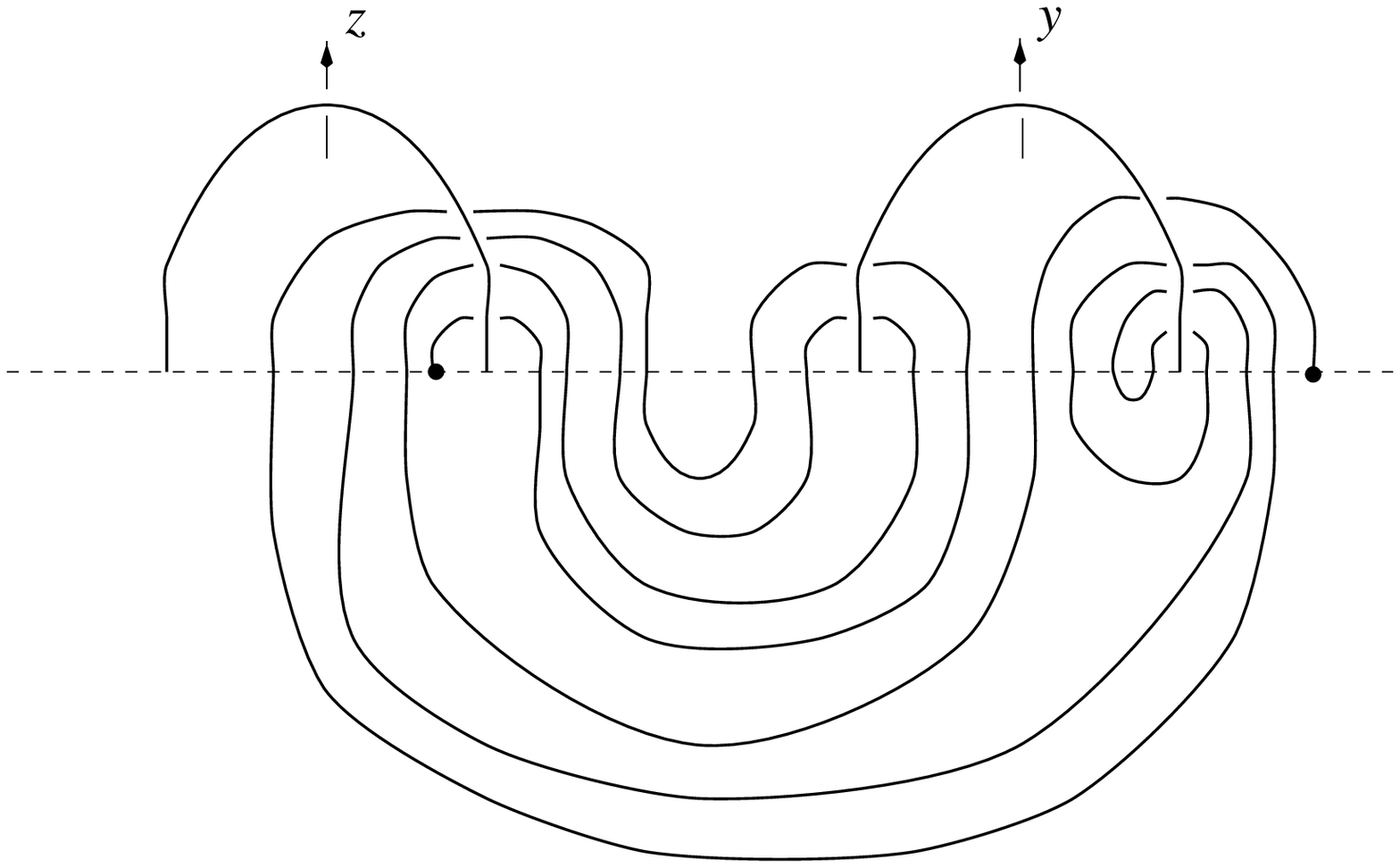} }
\medskip
\centerline {{\bf Figure 11.} {\msmall Both the endpoints
of the arc can be further hooked with the $z$-hook}}
\medskip

An example of an arc
where this exceptional 
case is realized is shown in Figure 11.
Here we have $i=1$. Notice that both endpoints of 
the arc shown here can be hooked with the $z$-hook.
Thus we can form $[y_1, y, y^{-1}, \ldots, y_j, \ldots]$,
where $y_1=z$, $y_j=z^{\pm 1}$ and the crossings
corresponding to $\{ y_j, \ y_j^{-1}\}$
occur on different endpoints of the $z$-hook.
\medskip
\p b) It follows immediately from part a) and 3.12(b).
\medskip
\p c) Let $d=[y_{1},\ \ldots,  
y_{j-1}]$ and suppose $y_{j}=y^{\pm 1}$,
for some free generator. Then
$[y_{1},\ \ldots,  
y_{j}]= [d, \ y^{\pm 1}]$. 
\smallskip
Let $\delta_1$
be the arc representing $[d, \ y^{\pm 1}]$.
A moment's thought will convince us that
in order for $\delta_1$ to be bad 
the following must be true: 
i) The arc $\delta_1^{*}$
of Lemma 3.10 must intersect the $y$-hook
precisely once; and
ii) crossings on the $y$-hook,
corresponding to some appearance of $y^{\pm 1}$
in $d$, must separate the crossings corresponding
to the canceling pair $\{ y_j, \ y_j^{-1}\}$.
In particular, $y^{\pm 1}$ must have appeared
in $d$ at least once.
Moreover it follows from 
Lemma 3.8
that, for any commutator $c$,
in order to be able to form
$[[c,\ y^{\pm 1}],\ y^{\pm 1}]$,
$[c,\ y^{\pm 1}]$ must be represented by a good arc.
Thus we 
may assume that $d$ satisfies the following:
at least one of the $y_{i}$'s 
is equal to $ y^{\pm 1}$; and  
$y_{j-1}\neq y^{\pm 1}$.

\p {\it Case 1:} Suppose that the arc
$\delta_1$
representing $[d, \ y_j]$
doesn't run through the hook of the band $B$;
in particular $\delta_1$ is embedded.
Then, any canceling pair in $d$ is admissible.
>From our assumption above,
the only remaining possibility
is when $\{ y_j, y_j^{-1} \}$ 
corresponds to the exceptional case
of part a). As already said in the proof of a),
in this case $\delta_1$ is good.

\p {\it Case 2:} Suppose that $\delta_1$ runs through
the hook of $B$. Then,

$$[d, \ y_j]= [e, \ y_r, \ \ldots, y_{j-1}, \ y_j],$$

\p where i) $\{ y_r, \ y_r^{-1} \}$
is the special canceling pair,
and
$\delta_1$ runs through the hook 
at this stage and ii)
$e= [y_1, \ldots, y_{r-1}]$ is a simple (quasi-)commutator.

Notice that all the canceling
pairs corresponding to $y_k$ with $k\neq r$
are admissible, by definition. 
Let $x_0$ be the free generator corresponding
to the band containing $\delta$.
By Lemma 3.8, we see
that
either $y_r=x_0^{\pm 1}$ or
the arc representing 
$[d, \ y_j]$ is embedded. In the later case,
it follows by Lemma 3.12 a), that
$\delta_1$ is embedded and the result follows as
in
Case 1.
If $y_r=x_0^{\pm 1}$ , then all the 
sub-arcs of
$\delta_1$ representing $e$
or $e^{-1}$ are embedded;
a  moment's thought will convince us
that are good arcs.
Thus if $j=r$, the conclusion of the lemma follows.
Suppose $j>r$.
By 3.12 a) we have $y_i\neq x_0^{\pm 1}$,
for all $i>r$. In particular, $y\neq x_0^{\pm 1}$.
Now the conclusion follows as in Case 1. \qed

\medskip

Before we are ready to state and prove the result 
about good arcs promised in the beginning of
3c, we need some more
notation and terminology.

Let
$c=[ y_1, \ y_2, \ \ldots, \  y_{n-1}, \  y_n]$
and let 
${\cal C}=\{ C_1, C_2, \ldots, C_{n-1}, C_n\} $ 
the sets of letters of Lemma 2.10. We will denote by $||\delta||$
the cardinality of the maximal subset of $\cal C$
that trivializes $\delta$ geometrically;
that is $\delta$ is $||\delta||$-nice (see Definition 3.4).
We will denote by
$s(c)$ the number of bad sets in $\cal C$. 
For a quasi-commutator ${\hat c}$, we will
define $s(\hat c)=s(c)$ where $c$
is the commutator underlying $\hat c$.

Finally, for $n\in \N$, let $t(n)$ be the quotient of the 
division of $n$ by four, and let $q(n)$
be the quotient of division by six.
\medskip
\p {\bf Lemma 3.15.}
{\sl Suppose that $S$, $B$, $\gamma$
and $\delta$ are as in the statement of Proposition 3.3,
and  that  $\delta_1$ is a good  subarc of $\delta$
representing a simple quasi-commutator $c_1$,  of length  $m+1$.

\p a) If $\delta_1$ is embedded then
 $\delta_1$ is an $t(m+1)$-nice arc.

\p b) If $\delta_1$ is non embedded
then  $\delta_1$ is  an  $q(m+1)$-nice arc.}
\medskip

\p {\it Proof.} \ {\sl a)} Inductively we will show that
$$||\delta_1 || \geq  m+1-s(c_1) \eqno (1)$$
Before we go on with the proof of (1),
let us show that it implies that $\delta_1$
is $t(m)$-nice.

For a
fixed free generator
$y$, let $w(y)$ be the number of appearances of $y$ in 
$c_1$ and let $s^{y}(c_1)$ be the number of bad
sets in $y$.
By Lemma 3.14 a), with one exception,
a set $C_i$
corresponding to $y$
can become bad only by a
successive appearance of $y$.
By Lemma 3.12, no letter can appear in $c_1$ more than
two successive times. 
A simple counting will convince us that
$${{w(y) \over  s^{y}(c_1)}}\geq {4 \over 3},$$
\p and that the maximum number of bad sets
in a
word is realized when each generator involved appears exactly four times, 
three of which are bad.
Thus we have $s(c_1) \leq m+1-t(m+1)$ and by (1)
we see that $||\delta_1 || \geq t(m+1)$, as desired.
 
We now begin the proof of (1), by induction on $m$.
For $m=1$, we know that 
all (embedded) arcs representing a  simple 2-commutator 
are nice and thus (1) is true. Assume that $m\geq 2$ and
(inductively) that for every good arc representing a
commutator of length $\leq m$, (1) is satisfied.
Now suppose that
$$c_1=[[c , \ z^{\pm 1}], \ y^{\pm 1}]$$
\p where $c $ is a simple commutator
of length  $m-1$, and $z, \ y$ are free generators. 

Let ${\bar {\delta}}_{1,2}$  (resp. ${\bar {\theta}}_{1,2, 3, 4}$)
denote the subarcs of $\delta_1$ representing
$[c^{\pm 1}, \ z^{\pm 1}]^{\pm 1}$ (resp. $c^{\pm 1}$).
\smallskip
\p {\it Case 1.} 
The arcs  ${\bar {\delta}}_{1, 2}$ are good.
By induction we have 
$$||{\bar {\delta}}_{1, 2} || \geq m-s(\bar c), \eqno (2)$$
\p where $\bar c = [c , \ z^{\pm 1}]$.
Since $\delta_1$ is good, a set of crossings that
trivializes ${\bar {\delta}}_{1, 2}$ can fail to
work for $\delta_1$ only if it becomes a
bad set in $c_1$. 
Moreover, the set of crossings
corresponding to the last canceling pair
$\{ y^{\pm 1}, \ y^{\mp 1} \}$ of $c_1$, also
trivializes $\delta _1$ geometrically. 
By 3.14a)
forming $c_1$ from  
$[c , \ z^{\pm 1}]$ can create at most two bad sets, each
involving
$y^{\pm 1}$. Thus we have 
$ s(\bar c)\leq s(c_1) \leq s(\bar c)+2$.
Combining all these with $(2)$, we obtain
$$||\delta_1||\geq ||{\bar {\delta}}_{1, 2} ||+1 -2\geq m+1-s (\bar c)-2
\geq m+1- s(c_1),$$
\p which completes the induction step in this case.
\smallskip
\p {\it Case 2.}
Suppose that ${\bar {\delta}}_{1, 2}$  are not good arcs.
 Let us use ${\bar {\theta}}$
to denote any of ${\bar {\theta}}_{1, 2, 3, 4}$.
Suppose that $c=[c_2,\  x]$, and thus 
$c_1=[[[c_2 , \ x], z^{\pm 1}] \ y^{\pm 1}]]$.
By Lemma 3.14(c),  ${\bar {\theta}}$ 
is a good arc and by induction
$$||{\bar \theta}|| \geq m-1-s(c). \eqno (3)$$

First suppose that $y\neq z^{\pm 1}$.
Because ${\bar {\delta}}_{1, 2}$ are not good
we can't claim that the pair
$\{ z^{\pm 1}, \ z^{\mp 1} \}$ trivializes
${\bar {\delta}}_{1, 2}$ geometrically;
however it will work for $\delta_1$.
Moreover,
the set of crossings
corresponding the last canceling pair
$\{ y^{\pm 1}, \ y^{\mp 1} \}$ 
of $c_1$, also
trivializes $\delta _1$ geometrically. 
Notice that the only sets of crossings
that work for $\bar \theta$ but could fail for
$\delta_1$ are these involving 
$z^{\pm 1}$ or $y^{\pm 1}$ that correspond to bad pairs in 
$\delta_1$. We see that
$$s(c)\leq s(c_1) \leq s(c)+4.$$
Combining all these with $(3)$, we obtain
$$||\delta_1||\geq ||{\bar \theta} ||+2-4 \geq m+1-s (c_1),$$
\p which completes the induction step in this case.

Now suppose that $y=z$. In this case we can see
 that
$s(c)\leq s(c_1) \leq s(c)+2$
and that at least one of the two last canceling pairs
$c_1$ will trivialize $\delta_1$ geometrically.
These together with (3) imply (1).
This finishes the proof of part {\sl a)} of our lemma.
\medskip
\p {\sl b)}  Let $x_0$ denote the free generator
of $\pi$ corresponding to it.
If $c_1$ doesn't involve $x_0$ at all, 
$\delta_1$ has to be
an embedded good arc and the conclusion follows
from part { \sl a)}. So we may suppose that
$\delta_1$ involves $x_0$.
Now the crossings that
correspond to  appearances of $x_0$
in $c_1$ may fail to trivialize
the arc geometrically.
\medskip
\centerline{\epsfysize=3in\epsfbox{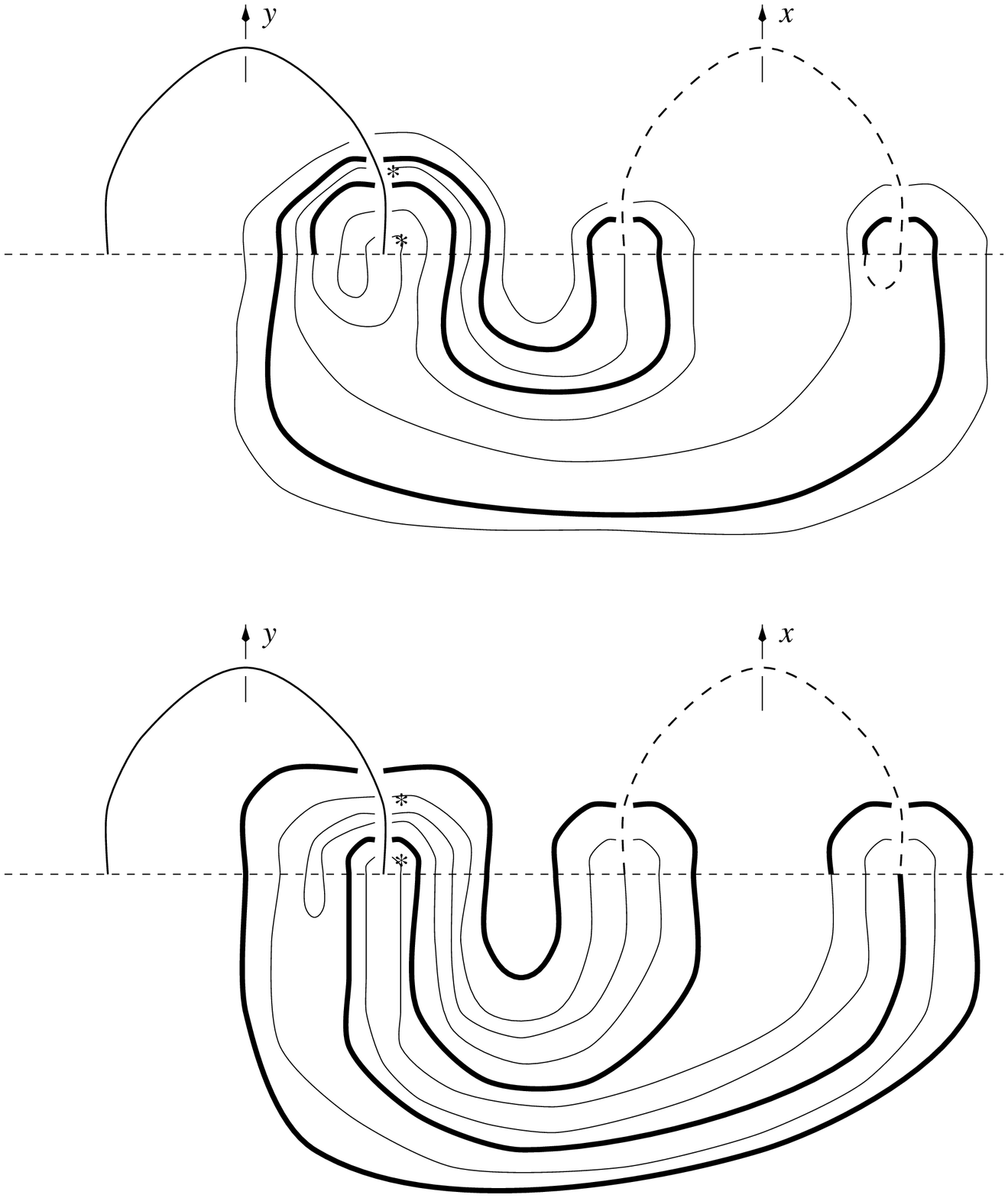} }
\medskip
\centerline {{\bf Figure 12.} {\msmall Simple commutators occupying
the entire band}}
\medskip

See, for example the arcs in
Figure 12;
in both cases 
crossings that realize 
the contributions of $x_0$
fail to trivialize the band.
As a result of this, we can only claim that
$$|| \delta_1|| \geq m+1- (s^{*}(x_0)+w(x_0)), \eqno (4)$$
\p where $w(x_0)$ denotes the number of appearances of 
$x$ in $c_1$ and
$s^{*}(x_0)$ is the number 
of the bad sets in generators different than $x_0$.

The proof of (4) is similar to that of (1) in part 
{\sl a)}.
Now by Lemmas 3.12 and 3.14 it follows that
the word $c_1$ that will realize the maximum
number of bad sets has the following parity:
$$[x_0, y_1, y_1, x_0, x_0,  \ldots, y_k, y_k,
x_0, x_0, y_1, y_1, \ldots, y_k, y_k],$$
\p where $y_1, \ldots, y_k$ are distinct
and $x_0\neq y_i$. Moreover,
three out of the four appearances of
each $y_i$ correspond to bad sets.
Now the conclusion follows. \qed
\medskip

\p {\bf Remark 3.15*.} From the proof of Lemma 3.15 we see that 
if the arc $\delta_1$
realizes the maximum number of bad sets then $q(m+1)=k$;
the number 
of distinct
generators, besides $x_0$, involved in 
$c_1$. In fact,  a careful reading of
the proof
of 3.15 will reveal
that
$$||\delta_1||\geq \cases {q(m+1), &if $m < 6k$\cr
\phantom{a} & \phantom{a} \cr
 k+\left[ \displaystyle { {{m-6k}\over 2} }\right], &if $m \geq 6k$.\cr }$$
Thus $||\delta_1||\geq q_{\delta_1}$,
and Proposition 3.3 follows
for bands representing
a simple (quasi-)commutator.
\medskip
\medskip
\p {\bf 3d. Conflict sets and products of good arcs}
\medskip
In this paragraph we study arcs 
that decompose into
products of good arcs. The main result of the paragraph
is Lemma 3.19, in which we show that an arc ${\tilde \delta}$
that is a product of good arcs is $q_{\tilde \delta}$-nice.

First, we need some more notation and terminology.
Let $S$, $B$, $\gamma$ and $\delta$
be as in the statement of Proposition 3.3,
and let ${\tilde \delta}$ be a subarc of $\gamma$ 
representing a word $W_1$ in $\pi^{(m+1)}$. 
Suppose that $W_1={\hat c}_1\ldots {\hat c}_s$,
is a product of quasi-commutators
represented by arcs $\{ \delta_1, \ldots, \delta_s\}$,
respectively.
Also, for $k=1,\ldots s$,
let ${\cal C}_k=\{C_{1}^k,\ldots, C_{(m+1)}^k\}$
be the sets of crossings of Lemma 2.10 for
$\delta_i$.
Let $C \in 2^{{\cal C}_k}$; by assumption 
the set of letters in $C$
trivialize $W_1$ algebraically.
For a
proper subset
${\cal D}\subset
\{ \delta_1, \ldots, \delta_s \}$ we will
use $C \cap {\cal D}$ (resp. $C \cap {\bar {\cal D}}
$)
to denote the crossings in $C$ that lie on arcs
in $\cal D$ (resp. in ${\bar {\cal D}}$
).
Here, ${\bar {\cal D}}$ denotes the complement
of $\cal D$ in the set
$\{ \delta_1, \ldots, \delta_s \}$ .

For every free generator, say $y$, we may have
crossings on the $y$-hook, realizing
letters in the word $W_1$, that trivialize geometrically some
of the subarcs $\delta_i$ but fail to trivialize $\tilde \delta$.
To illustrate how this can happen, consider the arcs $\delta_1$
and $\delta_2$. Let $C^{1}$ and $C^{2}$ be sets of crossings,
on the $y$-hook, along 
$\delta_1$
and $\delta_2$, respectively. 
Suppose that $C^{i}$ trivializes
$\delta_i$ geometrically (i.e. it
is
a
good set of crossings). Suppose,
moreover, that there are crossings on $\delta_2$
corresponding to a 
canceling pair $\{ y^{\pm 1}, \ y^{\mp 1}\}$
that doesn't 
belong in  $C^{2}$, and such that
they are separated by crossings in $C^{1}$. 
Then $C^{1}\cup C^{2}$ may not trivialize
$\delta_1 \cup \delta_2$.
With the situation described above in mind,
we give the following definition.
\medskip
\p {\bf Definition 3.16.} {\sl  A set $C$ of crossings
on $\tilde \delta$ is called {\sl a conflict set}
iff
i) the letters in $C$ trivialize $W_1$ algebraically; 
ii) switching the crossings in $C$ 
doesn't trivialize $\tilde \delta$ geometrically;
and
iii) there exists a proper subset
${\cal D}_C \subset
\{ \delta_1, \ldots, \delta_s\}$ such that
$C \cap {\cal D}_C$ trivializes geometrically
the union of arcs in ${\cal D}_C$ and 
$C \cap {\bar {\cal D}}_C$ trivializes geometrically
the union of arcs in  ${\bar {\cal D}}_C$.}
\medskip
Before we are ready to treat compositions of good arcs
we need two auxiliary lemmas that
describe the situations that create
conflict sets.

\medskip
\p {\bf Lemma 3.17.} {\sl For $k=1,\ldots s$,
let ${\hat c}_k$
be simple (quasi-)commutator represented by
an arc $\delta_k$, 
and let $\{ C_{1}^k,\ldots, C_{(m+1)}^k \}$
be sets of crossings as above.
Moreover
let $[y_{1}^k,\  y_{2}^k
,\  y_{3}^k, \ldots,  y_{m+1}^k]$ be
the underlying commutator of
${\hat c}_k$. Suppose that 
for some $i=1,\ldots, m+1$,
$y_i^k=y^{\pm 1}$ for some free generator $y$,
and that
$C_i= \cup_{k=1}^s C_i^k,$
is a conflict set.
Let ${\cal D}_{{C_i}}$ be as in Definition 3.16.
Then, there exist arcs $\delta_t \in {\cal D}_{{C_i}}$
and $\delta_r \in {\bar {\cal D}_{{C_i}}}$
such that the crossings on the $y$-hook
corresponding to canceling pairs in $y_j^r$ on $\delta_r$,
are separating by crossings corresponding
to $y_i^t$ on $\delta_t$. Here $j\neq i$.}
\medskip
{\it Proof.} It follows from
Definitions 2.3, 3.4 and 3.16. \qed
\medskip
With the notation as in Lemma 3.17,
the set $C_j^r$ will be called
{\sl a conflict partner} of $C_i^t$.
\medskip
\p {\bf Lemma 3.18.} {\sl 
Let 
$W_{1, 2}=[y_{1}^{1,2},\  y_{2}^{1,2}
,\  y_{3}^{1,2}, \ldots,  y_{m+1}^{1,2}]$ 
be the underlying commutators
of quasi-commutators represented by 
arcs
$\delta_{1,2}$. Let $C_{i}^1$ and $C_{j}^2$ 
be sets of letters in $W_1$ and $W_2$
respectively,
corresponding to the same free generator $y$.
Suppose that $C_{j}^2$ is a conflict partner
of
$C_{i}^1$.
Then, with at most one exception,

\p (1) either $j=i+1$
(resp. $j=i-1$) and $y_k^1=y_k^2$ for $k<i$
(resp. for $k< i-1$);
or

\p (2) the sets of free generators appearing
in $\{ y_{1}^1, \ldots, y_{i}^1 \}$
(resp. $\{ y_{1}^1, \ldots, y_{j}^1 \}$)
and $\{ y_{i+1}^2, \ldots, y_{j-1}^2 \}$
(resp. $\{ y_{j+1}^2, \ldots, y_{i-1}^2 \}$), 
are disjoint.}
\medskip
\p {\it Proof.} By 3.17, there must be crossings on
the $y$-hook corresponding to canceling pairs
in $C_j^2$, that are separated by crossings in $C_i^1$.
Let $C_1$ denote the canceling pair
corresponding to
$y_i^1$ in $[y_{1}^{1},\  y_{2}^{1}, \ldots,  y_{i}^{1}]$.
and
let $C_2$ denote the canceling pair
corresponding to $y_i^2$ (resp. $y_j^2$)
in $[y_{1}^{2},\  y_{2}^{2}, \ldots,  y_{i}^{2}]$
(resp. $[y_{1}^{2},\  y_{2}^{2}, \ldots,  y_{j}^{2}]$)
if $j>i$ (resp. $j<i$).
By Lemma 3.12({\sl a}), $C_{1, 2}$
are of type (I) or (II). Let $D_{1, 2}$
be the finite disc corresponding to
the canceling pair $C_{1,2}$,
in $W_{1, 2}$, respectively.
Up to symmetry there are three cases to consider:
i) Both $C_{1,2}$ are of type (I);
ii) both $C_{1,2}$ are of type (II); and
iii) one of them is of type (I) and 
the other of type (II).
In each case the result will follow
using Digression 3.9 to study 
the components of $D_{1, 2} \cap {\bar \delta_{1,2}}$,
where ${\bar \delta_{1,2}}$ denotes
the complement in ${\delta_{1,2}}$ of
the parallel arcs corresponding
to $C_{1, 2}$, respectively.
The exceptional case will occur
when the canceling pair $C_{1}$ is of type
(I) and crossings in it are a
separated by a type (II) canceling pair on
$\delta_2$.
The details are similar to these in the proof
part {\sl a)}
of Lemma 3.14  \qed
\medskip
To continue, recall the quantity
$q_{\tilde \delta}$ defined before
the statement of Proposition 3.3.
\medskip
\p {\bf Lemma 3.19.}\ (Products
of good arcs)  {\sl
 Let $S$, $B$, $\gamma$
and $\delta$ be as in the statement of Proposition 3.3,
and let $W=c_1\ldots c_r$ be
a word 
expressing $\delta^{+}$ as a product
of simple quasi-commutators.  Suppose that
${ \tilde \delta}$ is a subarc
of $\delta$, representing 
a subword of simple quasi-commutators 
$W_1={\hat c}_1\ldots {\hat c}_s$,
each of which is represented by a {\it good} arc.
Then  ${ \tilde \delta}$ is $q_{\tilde \delta }$-nice. In particular
if $W$ is a product of 
simple quasi-commutators
represented by good arcs,
Proposition 3.3 is true for ${ \delta}$.}
\medskip
\p {\it Proof.} 
If $s=1$ the conclusion follows from Lemma 3.15 and Remark 
$3.15^{*}$.
Assume that $s>1$. Let $\delta_1,\ \ldots, \ \delta_s$ 
be arcs representing  ${\hat c}_1\ldots {\hat c}_s$, respectively.

In general, we may have conflict sets of crossings between the $\delta_j$'s.
Since conflict sets occur between commutators
that have common letters, we must partition 
the set $\{{\hat c}_1, \ldots, {\hat c}_s\}$ into groups
involving disjoint sets of generators and work with each 
group individually. The maximum
number of conflicts will occur when
all the ${\hat c}_i$'s belong in the same group.
Since conflict sets
are in one to one correspondence
with proper subsets of $\{ \delta_1, \ldots, \delta_s\}$,
the maximum number
of conflict sets, for a
fixed generator $y$, is $2^{s}-2$. 

Let $x_0$ be the generator 
corresponding to the hook of $B$.
>From the proof of Lemma 3.15, Remark $3.15^{*}$
and by
Lemma 3.18 
we can see that a word $W$,
in which there are $k$ distinct generators besides $x_0^{\pm 1}$,
will realize the maximum number of bad sets of crossings
on the individual $\delta_i$'s and the maximum
number of conflict sets,
if the following are true:

\p i) The length $m+1$ is equal to $6k+r+2k(2^s-2)$,
where $r>2$;

\p ii) each of the arcs $\delta_i$ realizes
the maximum number of bad sets and the maximum number of appearances 
of $x_0^{\pm 1}$ (i.e. $5k+{\displaystyle {r\over 2}}$)
and there are $k(2^s-2)$ conflict sets
between the $\delta_i$'s. Moreover,
each pair of conflict partners in $W$
correspond either to the exceptional
case or in case (1) of Lemma 3.18.

We claim, however, that there will be
$k+{\displaystyle {\left[ r\over 2\right] }}+k(s-2)$ sets of crossings that trivialize
$\tilde \delta$ geometrically. From these $k+{\displaystyle {\left[
r\over 2\right] }}$ come from
good sets on the $\delta_i$'s.
The rest $ks-2k$ are obtained as follows:
For a fixed $y\neq x_0^{\pm 1}$,
the crossings in the conflict sets
involving $y^{\pm 1}$
and in their conflict partners 
can be partitioned into
$s-2$ disjoint sets that 
satisfy the definition of $(s-3)$-triviality.
To see that, create an $s\times (2^s-2)$
matrix, say $A$, such that the $(i,\ j)$
entry in $A$ is the $j$-th appearance
of $y$ in ${\hat c}_i$. The columns of $A$
are in one to one correspondence
with the conflict sets $\{C_i\}$, in $y$. By 3.18,
and 3.12 there are at most $2s$ ``exceptional"
conflict partners shared among the $C_i$'s. Other 
than that, the conflict partners of a
column $C_i$ will lie in exactly one of the  
adjacent columns.
For $s\geq 4$ we have $2^s-2\geq 4s$
and thus $A$ has at least $2s$
columns that can only conflict with
an adjacent column; these will give
$s>s-2$ sets as claimed above.
For $s=2, 3$ the conclusion is trivial.

Now from i) we see that
${\displaystyle {ks-2k>{\rm log}_2({{m+1-4k-r}\over {4}})}}$.
Thus, $${r\over 2}+k(s-2) > 
{\rm log}_2({{m+1-4k}\over {4}})
> 
{\rm log}_2({{m+1-6k}\over {6}}),$$
\p and
the claim in the statement of the lemma follows.\qed

\medskip
\medskip

\p {\bf 3e. The reduction to nice arcs}
\medskip
Let $S$, $B$, $\gamma$
and $\delta$ be as in the statement of
Proposition 3.3.
Our goal in this paragraph is to finish the proof
of 3.3. We begin with the following lemma, which relates
the $q_\mu$-niceness of a subarc 
$\mu\subset\tilde\delta\subset\delta$ to the 
$q_{\tilde\delta}$-niceness of $\tilde\delta$.
\medskip
\p {\bf Lemma 
3.20.}\ (Products
of nice arcs)  {\sl
Let $S$, $B$, $\gamma$
and $\delta$ be as in the statement of Proposition 3.3.
Let
${ \tilde \delta}$ be a subarc
of $\delta$, representing 
a subword $W_1={\hat c}_1\ldots {\hat c}_s$,
where ${\hat c_i}$ is a product
of simple quasi-commutators represented by
an arc $\theta_i$. Suppose that $\theta_i$
is $q_{\theta_i}$-nice, for $i=1, \ldots, s$.
Then
${ \tilde \delta}$ is $q_{\tilde \delta }$-nice.}
\medskip
\p {\it Proof.} Once again we can have 
sets of crossings on $\tilde \delta$
that trivialize $W_1$,
and trivialize a subset of $\{ \theta_1, \ldots, \theta_s \}$
geometrically but fail to trivialize
$\tilde \delta$.
For $i=1, \ldots, s$,
let $m_i$ denote the number of simple
quasi-commutators in ${\hat c}_i$,
and let ${\cal D}_i$
denote the set of subarcs of $\tilde \delta$ representing them.
We notice that the maximum number of conflict
sets that we can have in $W_1$, is $k[2^{(m_1+\ldots m_s)}-2]$ where
$k$ is the number of distinct generators, different
than $x_0$, appearing in $W$. 
Now we may proceed as in the proof of Lemma 3.19. \qed
\medskip 
To
continue recall the notion of a quasi-nice arc (Definition 3.4).
Our last lemma in this section shows that
the notions of quasi-nice and $q_{\tilde \delta}$-nice
are equivalent.
\medskip
\p {\bf Lemma 3.21.} {\sl A quasi-nice subarc 
${\tilde \delta} \subset \delta$ that represents
a product $W=c_1\ldots c_r$ of quasi-commutators,
is $q_{\tilde \delta}$-nice.}
\medskip
{\it Proof.} Let $\delta_1,\ldots, \delta_r$ be the arc
representing $c_1, \ldots, c_r$, respectively.
Let ${\cal D}_g$ (resp. ${\cal D}_b$) denote the
set of all good (resp. not good) arcs in $\{ \delta_1, \ldots, \delta_r \}$.
Also let $n_g$ (resp. $n_b$)
denote the cardinality of 
${\cal D}_g$ (resp. ${\cal D}_b$).
If $n_b=0$,
the conclusion follows from
Lemma 3.19. 
Otherwise let  $\mu \in {\cal D}_b$,
be the first of the $\delta_i$'s 
not represented by a good arc.
Suppose it represents $c_{\mu}=
[c^{\pm 1}, \
y^{\mp 1}]$, where  $c^{\pm 1}$
is a simple quasi-commutator of length $m$, and $y$ a
free generator. Let ${\mu}^{1,\ 2}$ be the subarcs
of $\mu$ representing $c^{\pm 1}$. 
Since  $\mu$ is not good, $y$
must have appeared in $c$; thus the numbers
of distinct gererators in the words representing
$\mu$ and ${\mu}^{1,\ 2}$ are the same.
We va see that
$q_{{\mu}^1}=q_{{\mu}^{2}}=q_{{\mu}}$.
By 3.14, 
${\mu}^{1,\ 2}$ are good arcs and by 3.15 they are
$q_{{\mu}}$-nice. Let $\bar \mu= {\tilde \delta} \setminus \mu$,
and let $\bar \mu_{1,2}$ denote the two
components of $\bar \mu$.
By induction and 3.20, $\bar \mu_{i}$ is
$q_{{\bar \mu_i}}$-nice. Since $\mu$ is not good,
one of its endpoints lies inside the $y$-hook
and the other outside. Moreover the arc $\mu^{*}$
of Lemma 3.10, separates crossings corresponding
to canceling pairs on the $y$-hook.
Now a moment's thought will convince us that
at least one of ${\bar \mu_{1,2}}$
must have crossings on the $y$-hook.
A set of crossings that trivializes geometrically
$\theta_1= {\mu}^{1}\cup {\mu^{2}}$ and $\theta_2=
{\bar \mu_{1}} \cup {\bar \mu_{2}}$ 
will fail to trivialize
$\tilde \delta$ only if there are conflict
sets between $\theta_1$ and $\theta_2$.
A counting argument 
shows that the maximum number 
of conflict sets that can be on $\tilde \delta$ is
$k(2^r-2)$, where $k$ is the number of distinct
generators, different than $x_0$, in $W$. Now the
conclusion follows as in the proof of 3.20. \qed

\medskip
\p {\it Proof of Proposition 3.3.} It follows immediately
from 3.21 and the 
fact thet the arc $\delta$ in
the statement of 3.3 is
quasi-nice;
see discussion at the end of 3b. \qed

\medskip
\p {\bf Remark 3.22.} a) Theorem 3.2 is not true if we don't
impose any restrictions on the surface $S$ of Definition 3.1.
For example let $K$ be a positive knot
set $\pi_K= \pi_1(S^3\setminus K)$
and let $D_K$ denote the untwisted Whitehead double of $K$. 
Let $S$ be the
standard genus one Seifert surface for $D(K)$. 
Since $\pi_K^{(n)}=\pi_K^{(2)}$ for any $n \geq 2$, we see that
$S$ has a half basis realized by a curve that  if
pushed in the complement of $S$ lies in
$\pi_K^{(n)}$,
for all $n \geq 2$. On the other hand,
$D_K$ doesn't have all its Vassiliev invariants trivial since
 it has non-trivial
2-variable Jones polynomial (see for example [Ru]).
\medskip
b) The methods of this section were applied in [K-L]
to obtain relations between finite type and Milnor invariants
of links that are realized as plat-closures of pure braids.
\bigskip
\centerline {  4. {{\B \uppercase 
{Vassiliev invariants as obstructions to  }}}{\it n-}{\B \uppercase {
sliceness}}}
\medskip
Let $S$ be a {\it regular} Seifert surface.
In this section we show
that the Vassiliev invariants
 of the knot ${K=\partial S}$
are null-concordance obstructions
of links that can be derived from {\it regular} spines of $S$.

Let $B^4$ denote the 4-ball, and let
$L=L_1\cup \ldots \cup L_m$
be an $m$-component link in $S^3=\partial B^4$.
We begin by recalling, from [Co] or [O],
the definition of {\it n-slice}
(or $n$-null-cobordant).

\medskip
\p {\bf Definition 4.1.} ( [Co], [O])\ {\sl
The link $L$ is called {\it n-slice}
if there exist $m$ disjoint
connected surfaces $V=V_1\cup \ldots \cup V_m \subset B^4$,
with $\partial V_i=L_i$ and such that
the following is true:
For some trivialization of a tubular neighborhood
of $V$ in $B^4$, which restricts to the standard
trivialization of  a tubular neighborhood
of $L$ in $S^3=\partial B^4$, the composition
$$\pi_1(V_i) \longrightarrow \pi_V \longrightarrow
{ \pi_V}/{ \pi_V}^{(n)}$$
\p is trivial, for all $i$.
Here, $\pi_V=\pi_1(B^4 \setminus V)$.}
\medskip

As was shown in [Co] and [O], $n$-slice
is the geometric notion which is equivalent
to the vanishing of Milnor's invariants.
More precisely we have:

\medskip
\p {\bf Theorem 4.2.}( [Co], [O])\ {\sl
A link $L=L_1\cup \ldots \cup L_m$ is $n$-slice 
for all $n\in \N$ if and only if
all its Milnor invariants vanish.
That is, the homomorphism
$F \longrightarrow  \pi_{L}$,
given by any
choice of meridians of $L$,
induces isomorphisms
$$F/ F^{(k)}\longrightarrow { \pi_L}/{ \pi_L}^{(k)}$$
for all  $k\leq n+1$ and all $n\in \N$. Here $F$
is the free
group of rank $m$, and $\pi_L=\pi_1(S^3 \setminus L)$.}
\medskip

To continue, let $S$ be a { regular} Seifert surface of a knot $K$,
and let $\Sigma$ be a regular spine of $S$,
consisting of curves 
$\gamma_1, \beta_1, \ldots, \gamma_g, \beta_g$.
Let $${\bar \epsilon}:=(\epsilon_1, \ldots, \epsilon_{2g})$$
\p where each $\epsilon_i$ is equal
to $+$ or $-$. Finally, let
$$L^{{\bar {\epsilon}}}=: 
\gamma_1^{\epsilon_1}\cup \beta_1^{\epsilon_2}\cup \ldots
\cup \beta_g^{\epsilon_{2g}}$$
\p Clearly, $L^{{\bar {\epsilon}}}$ is a link with $2g$-components.
Let $\pi=\pi_1(S^3\setminus S)$ and let ${ \pi_{L^{{\bar \epsilon }} }}=
\pi_1(S^3 \setminus L^{{\bar \epsilon}})$. By assumption,
$\pi$ is is the free
group of rank $2g$. 
We will say that the spine $\Sigma$ is {\it admissible}
if the following is true: For some ${\bar \epsilon}$ as above,
the longitudes of $L^{\bar \epsilon}$
lie in $\pi_{L^{{\bar \epsilon}} }^{(n)}$ (for some $n \in \N$)
if and only if $\gamma_i^{\epsilon_i}, \ \beta_i^{\epsilon_{i+1}} \in 
\pi^{(n)}$. The link $L^{\bar \epsilon}$
will be called an {\it admissible spine link} for $S$.
The theorem below shows the Vassiliev invariants
of the boundary knot are obstructions
of null-cobordance for the links $L^{{\bar \epsilon}}$.

\medskip
\p {\bf Theorem 4.3.} {\sl Let $S$ be a { regular} Seifert surface of a knot $K$
and suppose that $L^{{\bar \epsilon}}$ is an admissible spine link
for $S$. If the Milnor
invariants of length $\leq n+1$ vanish for
$L^{{\bar \epsilon}}$ 
then all the Vassiliev invariants of orders $\leq l(n)$
vanish for $K$.}  
\medskip
\p {\it Proof.} It follows from the fact that if the Milnor
invariant of length $\leq n+1$ vanish for
$L^{{\bar \epsilon}}$ then its longitudes 
lie in $\pi_{L^{{\bar \epsilon}} }^{(n+1)}$ ([M2]), our discussion
above, and Theorem 3.2. \qed

It was conjectured (see [Co]) that a link  $L$
is $n$-slice for some $n\in N$
if and only if Milnor's invariants 
with length less or equal to $2n$ vanish.
For the proof of one direction of this conjecture
(namely, that $n$-slice implies the vanishing
of invariants of length $\leq 2n$) see [L2].
The proof of the other direction was given
by Igusa and Orr ([I-O]). Thus we obtain.
\medskip
\p {\bf Corollary 4.4.}
{\sl If $L^{{\bar \epsilon}}$ is $n$-slice
then the Vassiliev invariants of orders $\leq l(2n-1)$
vanish for $K$.}
\medskip

\p {\bf Remark 4.5.} It is known (see [N-S] for relevant discussion)
that for every $n\in \N$,
and every non-slice knot $K$ one can find slice (null-concordant) knots
whose Vassiliev invariants  of orders $\leq n$
are equal to these of $K$.
It is interesting to compare Theorem 4.3 (or Corollary 4.4)
with this fact. 
\bigskip
\centerline {  5. {{\B \uppercase 
{ More types of }}}{\it n-}{\B \uppercase {trivial knots}}}
\medskip
Using the calculation of the Alexander polynomial via Seifert matrices
one can see that the Alexander polynomial
of an $n$-hyperbolic knot is trivial.
On the other hand there are $n$-trivial knots
with non-trivial Alexander polynomial.
Our purpose in this section is to
generalize
the notion of
$n$-hyperbolic so that we can include 
knots with non-trivial Alexander polynomial.
The generalized notion is that of $n$-{\it unknotted}
(see Definition 5.5). Roughly speaking, a knot is called
{\it $n$--unknotted} if it bounds a Seifert surface
whose complement looks, modulo certain terms of its
fundamental group, ``simple" ({\it $n$--parabolic}).
We show that the existence of such a surface 
for a knot implies the vanishing of its Vassiliev invariants of orders
$\leq n$. See Theorem 5.6. To complement the notions 
of $n$-hyperbolic and $n$-parabolic we also
distinguish another special class of $n$--unknotted;
namely that of {n-elliptic}. We will have more to say about 
each of these classes of knots and their relation with each other
in $\S 6$.
\bigskip
\p {\bf 5a. Definitions}
\medskip

We begin by introducing some notation and terminology needed
to continue.
Let $K$ be a knot and let $S$ be a
genus $g$ {\it regular} Seifert surface of $K$, in disc-band form.
For a band  $A$  of $S$ let 
$\gamma_A$ and  $x_A$ denote the  core and the free generator
of $\pi=\pi_1(S^3\setminus S)$ corresponding to $A$, respectively.
Also, for a set of bands 
 $\cal A$ we will denote by
${\cal G}_{\cal A}$ (resp. ${\cal D}_{\cal A}$)  the normal subgroup
of $\pi$ generated by $\{x_A|\ A \in {\cal A}\}$
(resp. $\{x_A|\ A \in {\cal {CA}}\}$). Here ${\cal {CA}}$
denotes the complement of ${\cal  A}$ .

\medskip
We will say that two bands $A$, $A^{'}$
are {\sl
geometrically related} if there exists a sequence of bands
$$A =A_0\rightarrow A_1 \rightarrow
\ldots \rightarrow A_{r-1} \rightarrow A_r=A^{'},$$
\p such that, for $i=1, \ldots, r$,
$\gamma_{A_i}^{+}$ or $\gamma_{A_i}^{-}$
is geometrically linked to at
least one of $\gamma_{A_0}, \ldots, \gamma_{A_{i-1}}$.
Otherwise $A$ and $A^{'}$ will be called 
{\sl
geometrically unrelated}.
Moreover, we will say that
two sets of bands   $\cal A$ and $\cal B$ are
{\it geometrically unrelated} if
any two $A \in {\cal A}$ and
$B \in {\cal B}$
are geometrically unrelated.

 For a band $A$ whose core $\gamma_A$
represents an element in $\pi^{(m+1)}$
let $q_{\gamma_A}$ be the quantity defined 
in the beginning of $\S 3$.
\medskip
\p {\bf Definition 5.1.} {\sl Let $n\in \N$, with $n >1$.
A regular Seifert surface $S$ will be called
{\it n-elliptic}, iff there exist
two half bases $\cal A$ and $\cal B$, 
represented by circles in a regular spine $\Sigma$,
such that the following is true:
The sets  
$\cal A$ and $\cal B$ are 
geometrically unrelated and
for every $A \in {\cal A}$ 
there exist
$B \in {\cal B}$  and  $m_A, m_B \in \N$ with $q_{\gamma_A}+q_{\gamma_B}=n+1$
and such that we have
$$[\gamma_A^{\epsilon}] \in {\cal G}_{\cal A}^{(m_A+1)}$$
\p and
$$[\gamma_B^{- \epsilon}] \in {\cal G}_{\cal B}^{(m_B+1)} $$
\p  The boundary
of such a surface
will be called an {\it n-elliptic} knot.}
\medskip

\p {\bf Definition 5.2.} {\sl
Let $n\in \N$, with $n >1$.
A regular Seifert surface $S$ will be called
{\it n-parabolic}, iff there exist
two half bases $\cal A$ and $\cal B$, 
represented by circles in a regular spine $\Sigma$,
such that the following is true:

\p a) The sets  
$\cal A$ and $\cal B$ are 
geometrically unrelated.
Moreover, any two subsets
of $\cal A$ are geometrically unrelated.

\p b) Let
$s=s(K, S) \geq 1$ be the largest integer such
that
$K$ can be unknotted 
in $2^s-1$ distinct ways by untwisting along the bands
in $\cal A$. If $n>s$, for every $B \in {\cal B}$
we have that
$$[\gamma_B^{\epsilon}] \in {\cal G}_{\cal B}^{(m+1)}$$

\p where $m\in \N$, such that 
with $q_{\gamma_B}+s=n+1$.

\p The boundary
of such a surface
will be called an  {\it n-parabolic} knot.
The numbers $s$ and $t=n+1-s$ will be called the
simplicity and triviality of $K$, respectively.}
\medskip

\p {\bf Remark 5.3.} Our terminology for regular surfaces,
in Definitions 3.1, 5.1 and 5.2 has been motivated
by the forms of their Seifert matrices.
To explain this, let $S$ be a regular surface and let $M$
be the matrix $V+V^T$, where $V$ is a Seifert matrix
of $S$. We see that if $S$ is $n$-elliptic, then there exists
a regular basis of $H_1(S)$, with respect to which
$M$ is of the form
$$\left[\matrix{ J& 0 & ... & 0 \cr
.& . & ...&  . \cr
0& 0 & ... & J \cr  }\right]$$
where
$$J=\left[\matrix{ 0 & 1 \cr
1 & 0 \cr }\right].$$
Thus $M$ is an elliptic matrix.

If $S$ is $n$-hyperbolic then $M$ is congruent to a matrix of
the form
$$\left[\matrix{ 0 & P \cr
P & R \cr }\right]$$
\p for some $g \times g$ matrices
$P$, $R$.
Finally, if $S$ is $n$-parabolic
then $M$ is congruent to a matrix of
the form
$$\left[\matrix{ D_g & D \cr D & T \cr }\right]$$
where $D_g$ and
$D$ are diagonal $g \times g$ matrices
and $T=0$ if $n>>0$.

We have:

\medskip
\p {\bf Theorem 5.4.} {\sl An $n$-elliptic
or $n$-parabolic knot is
$n$-trivial.}
\medskip

\p {\it Proof.} In both
cases the proof will be by induction
on the genus of the Seifert surface,
with the desired properties.

\p a) First we assume that $K$ is $n$-elliptic. 

If the genus of the surface in Definition 5.1 is zero, the
conclusion follows trivially.
Otherwise, let us concentrate on a
pair of dual  bands $(A,\ B)$.
We can assume that $[\gamma_A^{\epsilon}]\neq 1$
and $[\gamma_B^{-\epsilon}]\neq 1$.

Since by assumption the sets $\cal A$ and $\cal B$
are {\it geometrically unlinked}
we may modify the projection of $S$ so that Proposition 3.3
applies to both bands, simultaneously.
To see this,
first we find a projection with respect to which
the requirements of Lemma 2.12 are satisfied for $A$.
Let $l$ be the horizontal
line associated
to this projection.
Then, we start sliding the overcrossings of
$B$ that are below $l$ (see Remark 2.5) till
$B$ is put in good position.
At the end of the procedure, our projection will 
have the following property: Every crossing
between $A$ (resp. $B$) and any other band
that occurs below $l$ will be an undercrossing
(resp. overcrossing).

These properties
together with Proposition 3.3
guarantee the following:
We have  $n+1$ sets of crossings
($C_1,\ldots, C_{q(m_A)}$ from $A$
and $C_{{{q(m_{A})}+1}},\ldots, C_{n+1}$ from $B$)
that unable us to write $K$, in ${{\cal V}_n}$,
as a linear combination of $2^{n+1}-1$
knots $\{K_C\}$,
with the following property: each $K_C$ is 
an $n$-elliptic knot that
bounds an  $n$-elliptic Seifert surface, of genus less than $g$.
Thus the induction step applies.

\p b) Now we suppose that $K$ is $n$-parabolic.

Let $S$ be an $n$-parabolic Seifert surface for $K$,
and let $g$, $s$ and $t$ be the genus, the simplicity
and the triviality of $S$, respectively.
By Definition 5.2 the twists along the bands in $\cal A$
provide us with $s$ sets of crossings satisfying the definition of
$s-1$-trivial. These together with the $t=n+1-s$
sets of crossings obtained by Proposition 3.3
for some $B \in {\cal B}$ will allow us to write 
$K$, in ${{\cal V}_n}$,
as a linear combination of $2^{n+1}-1$
knots $\{K_C\}$,
with the following property: each $K_C$ is either the trivial knot,
or an $n$-parabolic knot that
bounds an $n$-parabolic surface, of genus less than $g$. \qed

\bigskip
\p {\bf 5b. Mixed type knots}
\medskip
Let the notation be as in the beginning of $\S 5a$
and let $\beta \in {\cal G}_{\cal A}^{(m+1)}$,
for some set of bands $\cal A$
and $m\in \N$.
Write $\beta$ as a product,
$W_1 \ldots W_s$, of commutators
in ${\cal G}_{\cal A}^{(m+1)}$ and partition
the set $\{ W_1, \ldots W_s\}$
into disjoint sets, say ${\cal W}_1, \ldots {\cal W}_t$
such that: i) $k_1+\ldots + k_{t}=l$, where $k_j$
is the number of distinct generators involved
in ${\cal W}_j$ and ii) for $a\neq b$, the sets
of generators appearing
${\cal W}_a$ and ${\cal W}_b$ are disjoint.
Let $$k ={\rm min} \{k_1, \ldots, k_t\}.$$
We define $$q_{\beta}= 
\cases {q(n+1), &if $n < 6k$\cr
\phantom{a} & \phantom{a} \cr
k+ \left[ {\displaystyle {{\rm log}_2({
{{n+1-6k}\over 6})}}}\right], &if $n \geq 6k$.\cr }
$$
Before we can sate the main result in this section,
which is a generalization of Theorems 3.2 and 5.4,
we need the following:

\medskip
\p {\bf Definition 5.5.} {\sl Let $n\in \N$, with $n>1$.
We will say that the knot $K$ is
{\it $n$--unknotted}
iff it has a regular Seifert surface $S$, that
contains two half bases $\cal A$ and $\cal B$
represented by circles in a regular spine $\Sigma$,
such that the following are true:
For every $A_i\in {\cal A}$ there exists a
$B_i \in {\cal B}$ 
so that 
$$[\gamma_{A_i}^{\epsilon}] =  (x_{A_i})^{l_{A_i}} \chi_{A_i} \mu_{A_i}$$
and
$$[\gamma_{B_i}^{- \epsilon}] = \zeta_{B_i} \chi_{B_i} $$
where $\mu_{A_i} \in \pi$,
$l_{A_i} \in \Z$, and $\zeta_{B_i} \in {\cal G}_{\cal B}$.
Moreover, we have the following:

\p a) There exist integers $m_1, \ldots , m_g$,
such that $\mu_{A_i} \in \pi_{{\cal A}_{i-1}} ^{(m_i+1)}$,
and with $q_{\mu_{A_i}}=n+1$

for $i=1, \ldots, g$.
Here $\pi_{{\cal A}_{i-1}}$ and ${\cal A}_i$
are as in Definition 3.1.

\p b) There exist $m_{A_i}, m_{B_i} \in \N$  such that
$$\chi_{A_i}\in {\cal G}_{\cal A}^{(m_{A_i}+1)}$$
and
$$\chi_{B_i}\in {\cal G}_{\cal B}^{(m_{B_i}+1)},$$
\p with $q_{\chi_{A_i}}+q_{\chi_{B_i}}=n+1.$

\p c) For every $A_i$, $B_i$ as in a) , either $ \chi_{A_i}=\chi_{B_i}=1$
or  $\chi_{A_i}\neq 1$ and $\chi_{B_i}\neq 1$ and either 
$\zeta_{B_i}=(x_{A_i})^{l_{A_i}}=1$ or
$(x_{A_i})^{l_{A_i}}\neq 1$ and $\zeta_{B_i} \neq 1$.
Futhermore, 
$\chi_{A_i}\neq 1$ and $\chi_{B_i}\neq 1$ if and only if
$\zeta_{B_i}=\mu_{A_i}=(x_{A_i})^{l_{A_i}}=1$.

\p d) All the  $\zeta_{B_i}$'s lie in ${\cal G}_{\cal B}^{(m+1)}$,
for some $m\in \N$ with the following property:
Let $S_C$ be the Seifert surface obtained by modifying
$S$ along a set of band crossings
$C$, that trivializes all of the $ \mu_{A_i}$, $\mu_{B_i}$,
$\chi_{A_i}$ and $\chi_{B_i}$. Let $s$ be the simplicity
of  $\partial (S_C)$. Then, $q_{\zeta_{B_i}}+s=n+1$.}
\medskip
\medskip
\p {\bf Theorem  5.6.} {\sl Assume that $K$, is
$n$-unknotted, for some $n\in \N$.
Then $K$ is $n$-trivial.}
\medskip

\p {\it Proof.} Assume that $K$ and $S$ are as in
Definition 5.5 and let $g$ be the genus of $S$.
Suppose that 
$S$ is in disc-band form
and let $A_i \in {\cal A}$ and $B_i\in {\cal B}$ be a pair of
bands of $S$
whose cores correspond
to  $\gamma_{A_i}$ and $\gamma_{B_i}$, respectively.

We will say that the pair $(A_i, B_i)$ is an {\it essential pair}
if we have  $\chi_{A_i}\neq 1$ and $\chi_{B_i}\neq 1$.

Let $a$ be the number of $A_i$'s with
 $\mu_{A_i} \neq 1$
and let $e$ be the number of {\it essential} pairs of bands.
Finally, let $b$ be the number of the non-trivial $\zeta_{B_i}$'s.
We define the following
complexity function for $(K , S)$:
$$d=d(K, S)= (g, \ b,\ e, \ a).$$
We order the complexities lexicographically, and proceed by induction on $d$.
For $d=(0,\ \* , \  \* , \  \* )$ the conclusion is trivially true,
since $K$ is the unknot.

Observe that if $b=0$ then , by Definition 5.5, either
$K$ is the unknot or it can be reduced to an $n$-elliptic
of genus less than $g$; in both cases the conclusion
of the theorem follows. 
Suppose, now, that $b>0$.  Notice that if $e>0$ then
there is a pair of bands $A_i,\ B_i$ such that
$[\gamma_{A_i}^{\epsilon}] =   \chi_{A_i} $
and
$[\gamma_{B_i}^{- \epsilon}] = \chi_{B_i} $.
Then we can argue as in Theorem 5.4 a),
to show that $K$ is equivalent in ${\cal V}_n$
to a summation of $n$-uknotted  knots with strictly less complexity.
Thus we may assume that $e=0$. Part d) of Definition
5.5 and Proposition 3.3  guarantee the existence of $n+1$ sets of crossings
that can be used to show that $K$
is equal, in ${\cal V}_n$, to an alternating summation of knots
$K_C$  such that,  each $K_C$ bounds $n$-unknotted
surface that either has genus $<g$ or it contains a
band $A_i$ such that  $[\gamma_{A_i}^{\epsilon}]=
\mu_{A_i} \in \pi^{(m_i+1)} $, with
$q_{\mu_{A_i}}=n+1$. In both cases, it follows by Proposition 3.3
that $K$ is equal, in ${\cal V}_n$, to an alternating summation of knots
$K_C$  such that,  each $K_C$ 
is an $n$-unknotted with $d(K_C)< d(K)$. \qed

\bigskip
\smallskip
\centerline {6. {\B \uppercase {Discussion and questions}}}
\medskip

In $\S 3$ and $\S 5$ we studied various classes
of $n$-trivial knots.
The following Lemma, whose proof follows directly from
Definition 3.1, 5.1, 5.2 and 5.5, 
and the discussion following clarifies the relation
between these three types.
\medskip
\p {\bf Lemma 6.1.} {\sl a) Let $K$ be a $2n$-elliptic knot.
Then $K$ is $n$-hyperbolic.

\p b) Let $K$ be an $n$-parabolic knot with simplicity $s$,
and assume that $n > s+1$. Then $K$
is $(n-s-1)$-hyperbolic.

\p c) Let $K$ bne a $2n$-unknotted knot and let $s$
be the quantity defined in definition 5.5.
In $2n>s+1$ then $K$ is $(2n-s-1)$-unknotted. }
\medskip
Let ${\Delta_K}(t)$ denote the
Alexander polynomial of a knot $K$.
Using the representation of Alexander polynomial via Seifert matrices
(see, for example, [Ro]) one can see that
${\Delta_K}(t)=1$
for any $n$-hyperbolic knot $K$.
\bigskip
\centerline{\epsfysize=1in\epsfbox{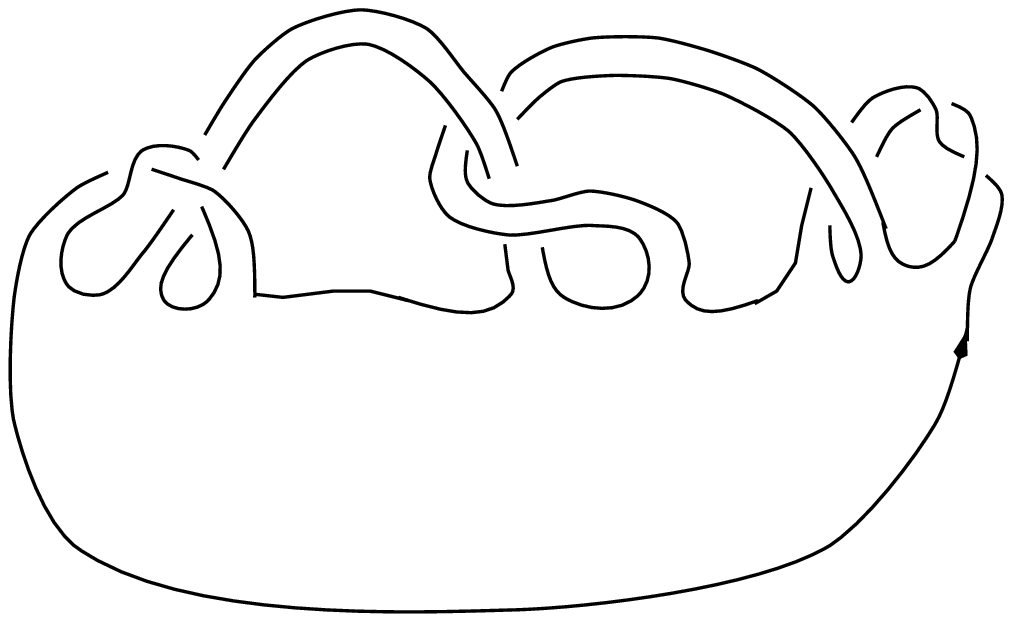} }
\medskip
\centerline{ { \bf Figure 13.} {\msmall A 2-parabolic knot}}
\bigskip
On the other hand
in Figure 13 we show a 2-parabolic knot,
with simplicity $s=2$,
that has nontrivial Alexander polynomial.
Moreover we have:
\medskip
\p {\bf Proposition 6.2.} {\sl For every $n \in \N$, there
is an $n$-trivial knot $K$, whose Alexander polynomial
is non trivial.}
\medskip

\p {\it Proof.} For a knot $K$, let $\Delta_K$
denote its Alexander polynomial.
By the Morton-Melvin expansion ([BN-G], [M-M])
for  the Alexander polynomial
we have
$${{p(h)}\over {\Delta_K(e^h)}}= \sum_{n} v_n(K) h^n$$
\p where $h$ is a variable and $p(h)={\displaystyle {e^{h\over 2}-
e^{-{h\over 2}}}\over {h}}$.

Moreover, each $v_n(K)$ is a {\it canonical} invariant of orders $n$,
determined by the colored 
Jones polynomial of $K$.
Now the result follows
from the fact that given a 
{\it canonical} (determined by its
weight system) Vassiliev invariant
$v$ of order $n+1$, there is an $n$-trivial knot $K$,
such that $v(K) \neq 0$ (see for example [N-S]). \qed
\medskip
Thus the Alexander polynomial of an $n$-trivial knot
obstructs to the knots being $m$-hyperbolic for any $m \in \N$.
Now let us start with a knot $K$ such that 
${\Delta_K}(t)=1$. This implies 
that $K$ is algebraically slice; that is 
for any Seifert surface $S$
we may find a basis of $H_1(S)$ with respect to which
the Seifert matrix of $S$ is of the form 
$$V= \left[\matrix{ 0 & A \cr
B & C\cr }\right]$$
\p for some $g \times g$ matrices
$A$, $B$ and  $C$. By a further change of basis
of $H_1(S)$, and using the condition 
${\Delta_K}(t)=1$, we may change $V$ 
to one of the forms that the linking matrix of
a 0-hyperbolic has. Thus, we see that 
$K$ is ``algebraically" 0-hyperbolic.
This suggests that the Alexander polynomial
might be the only obstruction to 
$m$-hyperbolicity of $n$-trivial knots. We conjecture:
\medskip
\p {\bf Conjecture 6.3.} {\sl 
A knot $K$ is $n$-trivial for all $n\in \N$
if and only if it is $n$-hyperbolic
for all $n\in \N$.}
\medskip
As already mentioned, the Alexander polynomial of a
knot is determined by linking numbers of 
``spine links" of Seifert surfaces. In particular, if
``enough" of these linking numbers vanish
then all the derivatives of the Alexander polynomial vanish.
Our Theorem 3.2 is a generalization of this fact.
It asserts that if ``enough" of the higher order linking
numbers (certain values Fox's higher order 
free derivatives)
of spines vanish then all the 
derivatives
of all the Jones type polynomials, up to certain order, vanish.
This in turn suggests that one might be able to
express certain Vassiliev invariants of a knot $K=\partial S$,
in terms of, appropriate values of, Fox derivatives
of push offs of regular spines.
For example, let $S$ be a
regular surface of $K$ and 
let $\Sigma$ be a regular spine of $S$,
consisting of curves 
$\gamma_1, \beta_1, \ldots, \gamma_g, \beta_g$.
Also, let ${\gamma_1^{\epsilon_1}, \beta_1^{\epsilon_2}, \ldots,
\beta_g^{\epsilon_{2g}}}$ be as in $\S 4$ and le
$W_1^{\epsilon_1}, W_1^{\epsilon_2}, \ldots,
W_{2g}^{\epsilon_{2g}}$ be words representing
them in the free group
$\pi_1(S^3 \setminus S)$.
One
can ask the following: Do the coefficients of the
Magnus expansions of the $W_i^{\epsilon_j}$'s
provide obstructions to the vanishing
of all Vassiliev invariants of $K$?
This question is very important for Conjecture 6.3;
in fact one needs to study ``modified" finite type invariants
of triples $(K,\ S,\ \Sigma)$, in terms 
Fox derivatives
of push offs of the spine $\Sigma$.
We will address this question in a subsequent preprint.

It would be interesting to know
to what extent our notion of $n$-unknotted 
captures the vanishing of Vassiliev's
invariants of bounded orders. We ask the following:
\medskip
\p {\bf Question 6.4.} {\sl 
Is it true that all the Vassiliev invariants of orders $\leq n$
of a knot $K$ vanish, if and only if it is $k$-unknotted
for some $k=k(n)$? Here $k=k(n)$ has to be an increasing 
function of $n$.}
\medskip
\p {\bf Remark 6.5.} A knot is called $n$-adjacent to the unknot if it is $n$-trivial
and  each set of crossings in the definition of $n$-triviality has cardinality one.
The first named author of this paper and N. Askitas
have proved  Conjecture 6.3 for $n$-adjacent knots. ([A-K]).

\bigskip
\bigskip
\bigskip

{\centerline{ \B REFERENCES}}
\medskip
\item {[A-K]} E. Kalfagianni and N. Askitas: {\it On $n$-adjacent knots},
preprint in preparation.

\item {[Bi]} J. S. Birman: {\it  New  points of  view
  in knot theory}, Bull. of AMS, {\bf 28}(1993), 253--287.

\item {[BN1]} D. Bar-Natan: {\it On the Vassiliev knot invariants}, 
Topology, {\bf 34}(1995), 423-472.

\item {[BN2]} D. Bar-Natan: {\it Vassiliev homotopy string link invariants},
J. of Knot Theory and Ramifications, {\bf 4}(1995), 13-32. 

\item {[BN-G]} D. Bar-Natan and S. Garoufalidis:
{\it On the Melvin-Morton-Rozansky conjecture}, Invent. Math., {\bf 125}(1996), 103--133.

\item {[Co]} T. D. Cochran: {\sl Derivatives of Links:
Milnor's Concordance Invariants and 
Massey's Products},
Memoirs AMS, {\bf 84} (1990).

\item {[C]} J. Conant: Private communication.

\item {[C1]} J. Conant:
{\it A knot bounding a grope of class n is n/2-trivial}, UCSD
preprint 1999.


\item {[Gu]} M. N. Gusarov:
{\it On $n$-equivalence of knots and invariants of finite degree},
in {\sl Topology of Manifolds and Varieties}, 173--192, 
Adv. Soviet Math., vol. {18}, AMS, 1994.

\item {[He]} J. Hempel: {\sl 3-manifolds},
Annals of Mathematics Studes, vol. 86, Princeton University Press, 1976.

\item {[H-M]} N. Habegger and G. Masbaum: {\it
The Kontsevich integral and Milnor's Invariants},
preprint 1997.

\item {[I-O]} K. Igusa and K. Orr:
{\it Links, pictures and the homology of
nilpotent
groups}, preprint, 1998.

\item {[KSM]} A. Karras, W. Magnus and D. Solitar: 
{\sl Combinatorial Group Theory}, Pure and Appl. Math.,
vol. XIII, Interscience, NY, 1966.

\item {[K-L]} E. Kalfagianni and X.-S. Lin: {\it
Milnor and finite type invariants of
 plat-closures}, Math. Reasearch Letters, to appear .

\item{[L1]} X.-S. Lin: {\it Power series expansions and invariants of links},
in {\sl Geometric Topology}, 184--202, AMS/IP Studies in Advanced Mathematics, 
vol. 2, 1996.

\item {[L2]} X.-S. Lin:
{\it Null k-cobordant links in $S^3$}, Comm. Math. Helv., {\bf 66}(1991), 333-339.

\item {[Mi1]} J. Milnor: {\it Link groups},
Annals of Math., {\bf 2} (1954), 145-154.

\item {[Mi2]} J. Milnor:
{\it Isotopy of links}, in {\sl Algebraic Geometry and Topology}, 280--386,
Princeton University Press, 1957.

\item {[M-M]} P. Melvin and H. Morton:
{\it The colored Jones function}, Comm. Math. Phys., {\bf 169}(1995), 501-520.

\item{[N-S]} K.-Y. Ng and T. Stanford: {\it On
Gussarov's groups of knots}, Math. Proc. Camb. Phil. Soc., to appear.

\item {[O]} K. Orr: {\it Homotopy invariants
of links}, Invent. Math., {\bf 95}(1989), 379-394.

\item {[Ro]} D. Rolfsen: {\it Knots and Links},
Publish or Perish, 1976.

\item{[Ru]} L. Rudolph: {\it A congruence between link polynomials}, 
Math. Proc. Camb. Phil. Soc.,
{\bf 107} (1990), 319-327.

\item{[Ya]} M. Yamamoto: {\it Knots and spatial embeddings of the complete 
graph on four vertices,}
Top. and Appl., {\bf 3}(1990), 291-298.
\bigskip

\p{\small {}}
{\Smaller {\small D}EPARTMENT OF {\small M}ATHEMATICS,
 {\small M}ICHIGAN {\small S}TATE
{\small U}NIVERSITY, {\small E}AST {\small L}ANSING, 
{\small MI}
{\small 48824}}

\p {\small { E--mail address:}}
{\vsmall kalfagia@math.msu.edu}
\medskip
\p{\Smaller {\small D}EPARTMENT OF {\small M}ATHEMATICS,
 {\small U}NIVERSITY OF
{\small C}ALIFORNIA, {\small R}IVERSIDE, 
{\small CA}
{\small 92521}}

\p{\small {E--mail address:}} {\vsmall  xl@math.ucr.edu}

\end

\input epsf
\end